\documentclass[12pt]{amsart}

\usepackage[dvipsnames]{xcolor}
\usepackage[utf8]{inputenc}
\usepackage{amsmath,amsthm,amssymb}
\usepackage{indentfirst}
\usepackage{url}
\usepackage{tikz-cd}
\usepackage[colorlinks]{hyperref}
\usepackage{colonequals}
\usepackage{graphicx}
\usepackage{enumerate}
\usepackage{soul}
\usepackage{palatino}
\usepackage{mathrsfs,mathtools}
\usepackage{bbm}
\usepackage[alphabetic]{amsrefs} 
\usepackage{ytableau}

\usepackage[margin=0.8in]{geometry}

\newtheorem{thm}[subsection]{Theorem}
\newtheorem{lem}[subsection]{Lemma}
\newtheorem{prop}[subsection]{Proposition}
\newtheorem{cor}[subsection]{Corollary}
{
\theoremstyle{definition}
\newtheorem{remx}[subsection]{Remark}

\newtheorem{defnx}[subsection]{Definition}
\newtheorem{examplex}[subsection]{Example}
}
\newenvironment{rem}
{\pushQED{\qed}\remx}
{\popQED\endremx}

\newenvironment{defn}
{\pushQED{\qed}\defnx}
{\popQED\enddefnx}
\newenvironment{example}
{\pushQED{\qed}\examplex}
{\popQED\endexamplex}

\newcommand{\NN}{\mathbb N}
\newcommand{\ZZ}{\mathbb Z}
\newcommand{\QQ}{\mathbb Q}
\newcommand{\RR}{\mathbb R}
\newcommand{\CC}{\mathbb C}
\newcommand{\PP}{\mathbb P}

\newcommand{\conv}{\mathrm{conv}}
\newcommand{\ext}{\mathrm{ext}}
\newcommand{\Ext}{\mathrm{Ext}}
\newcommand{\init}{\mathrm{in}}
\newcommand{\red}{\mathrm{red}}

\newcommand{\Cay}{\mathrm{Cay}}

\DeclareMathOperator{\Vol}{Vol}
\DeclareMathOperator{\covol}{covol}
\DeclareMathOperator{\codim}{codim}
\DeclareMathOperator{\Hilb}{Hilb}

\DeclareMathOperator{\Newton}{Newton}
\DeclareMathOperator{\supp}{supp}

\DeclareMathOperator{\bideg}{bideg}
\DeclareMathOperator{\mdeg}{mdeg}

\definecolor{darkred}{rgb}{0.7,0,0}
\newcommand{\demph}[1]{{\color{darkred}\emph{#1}}}

\newcommand{\indeg}{\mathrm{indeg}}
\newcommand{\outdeg}{\mathrm{outdeg}}

\begin{document}
\title{Generalized Alexander polynomials via Gr\"obner geometry}
\author{Linus Setiabrata}
\address{Linus Setiabrata, Department of Mathematics, Massachusetts Institute of Technology, Cambridge, MA 02139.}
\email{setia@mit.edu}
\date{\today}
\maketitle

\begin{abstract}
We study the Gr\"obner geometry of a toric ideal $I_{L(\mathcal A)}$ designed to geometrize the generalized Alexander polynomial $f_{\mathcal A}(t)$ of a flat vector arrangement $\mathcal A$. Our work implies that the coefficients of $f_{\mathcal A}(t)$ form a log-concave sequence whenever $\mathcal A$ is full-dimensional. Inspired by Gr\"obner geometry in Schubert calculus, we focus on graphical arrangements $\mathcal A = \mathcal A_G$. In this case, $I_{L(\mathcal A_G)}$ defines a subvariety $\widehat X_G$ of a product of projective spaces indexed by the vertices of $G$. We show that the multidegree $\mathfrak S_G(\mathbf z)$ of $\widehat X_G$ is the reciprocal of the integer point transform of the trimmed zonotope $Z_G^-$. In particular, $\widehat X_G$ is multiplicity-free, so the twisted $K$-polynomial $\mathfrak G_G(\mathbf z)$ has support equal to a generalized polymatroid. We show that the top-degree part of $\mathfrak G_G(\mathbf z)$ is the integer point transform of the interior of the zonotope $Z_G$.
\end{abstract}

\section{Introduction}

Let $\mathcal A = \{a_1, \dots, a_N\}$ be a full-dimensional vector arrangement in $\RR^n$, and assume that it is \emph{flat}, i.e., that there is a linear functional $h\colon\RR^n\to\RR$ so that $h(a_i)=1$ for all $i$. Li--Postnikov \cite{lp13} associated to $\mathcal A$ the \emph{generalized Alexander polynomial}
\begin{equation}
\label{eqn:alex}
f_{\mathcal A}(t)=\sum_{B\in\mathcal B(\mathcal A)} \Vol(\Pi_B)t^{\ext_\rho(B)},
\end{equation}
where $\mathcal B(\mathcal A)$ denotes the set of bases of the matroid defined by $\mathcal A$, $\Pi_B$ is the parallelepiped spanned by the vectors in $B$, and $\ext_\rho(B)$ is the external semi-activity of the basis $B$ with respect to a \emph{framing vector} $\rho\in\RR^N$. The polynomial $f_{\mathcal A}(t)$ is independent of $\rho$.

The polynomial $f_{\mathcal A}(t)$ is a rich invariant of the flat vector arrangement $\mathcal A$, and is closely related to many aspects of the polyhedral geometry of the zonotope $Z_{\mathcal A}$ spanned by $\mathcal A$. For example, they were originally studied \cite{lp13} in connection with the combinatorics of fine mixed subdivisions of $Z_{\mathcal A}$, and can be used to study the volume of $Z_{\mathcal A}$ together with the integer points of the trimmed zonotope $Z_{\mathcal A}^-$. This theory vastly extends Postnikov's study \cite{postnikov09} of triangulations and volumes of root polytopes $Q_G$ along with the integer points of trimmed $y$-generalized permutahedra $P_G^-$: in this language, $Q_G$ is $Z_{\mathcal A} \cap \{h(x) = 1\}$ and $P_G^-$ is $\{x \in Z_{\mathcal A}^-\colon h(x) \textup{ is minimized}\}$.

This paper grew out of our study of an algebraic variety which is designed to capture information about $f_{\mathcal A}(t)$. Our approach is inspired by Gr\"obner geometry in Schubert calculus \cite{km05}, although the tools we shall need turn out to predate \cite{km05}. We study a $\ZZ^2$-graded ideal $I_{L(\mathcal A)}$ defining a biprojective variety $X_{\mathcal A}$ which is equipped with a Gr\"obner degeneration $X_{\mathcal A}\rightsquigarrow X_0$. The limit $X_0$ is an arrangement of coordinate subspaces manifestly possessing combinatorics encoding the formula~\eqref{eqn:alex} expressing $f_{\mathcal A}(t)$ as a sum of monomials, one for each element of a combinatorially interesting set $\mathcal B(\mathcal A)$. It follows that $f_{\mathcal A}(t)$ represents the Chow class of $X_{\mathcal A}$. In this sense, $X_{\mathcal A}$ geometrizes $f_{\mathcal A}(t)$.

\subsection*{Graphical arrangements}
Our principal example is the graphical arrangement $\mathcal A = \mathcal A_G$ for a graph $G$. In this case, the ideal $I_{L(\mathcal A_G)}$ admits a multigrading by $\ZZ^{V(G)}$, and consequently defines a subvariety $\widehat X_G$ of a product of projective spaces indexed by the vertices of $G$. We study the multidegree $\mathfrak S_G(\mathbf z)$ of $\widehat X_G$, inspired by the combinatorial theory of Schubert polynomials.

This construction applies even when $\mathcal A = \mathcal A_G$ is not flat. However, when $G$ is bipartite and the edges of $G$ are oriented from the left part to the right part (so that $\mathcal A_G$ is flat), the $\ZZ^{V(G)}$-multigrading on $I_{L(\mathcal A)}$ refines the bigrading used to produce the biprojective variety $X_{\mathcal A}\subseteq \PP^N \times \PP^N$ discussed earlier, and it follows that a certain specialization of $\mathfrak S_G(\mathbf z)$ equals the generalized Alexander polynomial $f_{\mathcal A}(t)$.

As above, we have a Gr\"obner degeneration $\widehat X_G \rightsquigarrow \widehat X_0$, and we compute explicitly the multidegree $\mathfrak S_G(\mathbf z)$ by matching it with a description \cite[Theorem 1.4]{ms26} (implicit also in \cite[Theorem 6.5]{kmp25}) of the set of integer points of trimmed graphical zonotopes.

\begin{thm}
\label{thm:S_G}
Let $G$ be any connected directed graph. Then $\mathfrak S_G(\mathbf z)$ is the reciprocal of the integer point transform of a suitable translate of $Z_G^-$. More precisely,
\[
\mathfrak S_G(\mathbf z) = \mathbf z^\alpha \cdot \left(\sum_{\mathbf p \in Z_G^- \cap \ZZ^{V(G)}} \mathbf z^{-\mathbf p}\right),
\]
where $\alpha \colonequals \outdeg(G) - \mathbf 1_{V(G)} \in \ZZ^{V(G)}$ is a shift of the outdegree vector of $G$.
\end{thm}

We also study the twisted $K$-polynomial $\mathfrak G_G(\mathbf z)$, with inspiration from the combinatorial theory of Grothendieck polynomials. The Gr\"obner limit $\widehat X_0$ turns out to be reduced, and the Stanley--Reisner complex is a triangulated ball, so the twisted $K$-polynomial $\mathfrak G_G(\mathbf z)$ can be computed combinatorially \cite[Proposition 4.2]{kmy08} (cf.\ also \cite[Theorem 4.1]{km04}). One consequence is that the top-degree part $\mathfrak G_G^{\mathrm{top}}(\mathbf z)$ of $\mathfrak G_G(\mathbf z)$ is also a certain integer point transform (up to sign) -- as far as we can tell, for reasons highly unrelated to Theorem~\ref{thm:S_G}(!).

\begin{thm}
\label{thm:G_G-top}
Let $G$ be a connected directed graph with no bridges. Then
\[
\mathfrak G_G^{\mathrm{top}}(\mathbf z) = (-1)^{|V(G)| - 1} \mathbf z^{\outdeg(G)} \cdot \left(\sum_{\mathbf p \in \mathrm{relint}(Z_G)\cap \ZZ^{V(G)}} \mathbf z^{-\mathbf p}\right),
\]
where $\mathrm{relint}(Z_G)$ denotes the relative interior of the zonotope $Z_G$.
\end{thm}
If $e$ is a bridge of $G$, then $\mathfrak G_G(\mathbf z) = \mathfrak G_{G\setminus e}(\mathbf z)$ and $Z_G = Z_{G\setminus e} \times [0,\mathbf 1_{s(e)} - \mathbf 1_{t(e)}]$, so the bridgeless assumption on $G$ is essentially harmless. The coefficients in intermediate degree components of $\mathfrak G_G(\mathbf z)$ may be strictly larger than $1$.

The formula in Theorem~\ref{thm:S_G} implies that $\widehat X_G$ is \emph{multiplicity-free} in the sense of Brion \cite{brion03}. In the fortunate setting of multiplicity-freeness, work of Knutson \cite[Theorem 3]{knutson09} and Castillo--Cid-Ruiz--Mohammadi--Monta\~no (\cite[Theorem A]{ccmm22}) gives strong control over $\mathfrak G_G(\mathbf z)$.
\begin{thm}
\label{thm:G_G-extra}
The homogenization $\widetilde{\mathfrak G}_G(\mathbf z)$ of the twisted $K$-polynomial of $\widehat X_G$ has M-convex support, and in particular satisfies the saturated Newton polytope property. Furthermore, the coefficients of $\mathfrak G_G(\mathbf z)$ can be computed via M\"obius inversion:
\[
\mathfrak G_G(\mathbf z) = \sum_{\mathbf q\in P_G}\mu(\mathbf q) \mathbf z^{\mathbf q},
\]
where $P_G$ is the poset $\{\mathbf q \in \ZZ_{\geq 0}^{V(G)}\colon \mathbf q \leq \deg(G) - \mathbf 1\textup{ and } \mathbf z^{\mathbf q}\textup{ is divisible by some monomial in } \mathfrak S_G(\mathbf z)\}$ and $\mu$ is the M\"obius function of $P_G$.
\end{thm}

\subsection*{Connection with Alexander polynomials of links}
Although the generalized Alexander polynomial $f_{\mathcal A}(t)$ appears naturally from the theory of oriented matroids and of zonotopes, there is also a surprising connection to knot theory introduced by K\'alm\'an--M\'esz\'aros--Postnikov \cite{kmp25}. When $\mathcal A = \mathcal A_G$ is the graphical arrangement for a planar bipartite graph $G$, the polynomial $f_{\mathcal A}(t)$ is equal to the Alexander polynomial $\Delta_{L_G}(-t)$ (with the sign substitution so that coefficients are positive) of the special alternating link $L_G$ associated to $G$.

Theorem~\ref{thm:log-concave} gives a new, independent proof of the log-concavity of the coefficients of the Alexander polynomial of special alternating links (\cite{hmv24a,hmv24b}). It would be interesting to find a direct link between the variety $\widehat X_G$ and the link $L_G$, and in particular whether the variety $X_{\mathcal A}$ encodes richer invariants of $L_G$. (We caution the reader that the graph $G$ itself is not an invariant of $L_G$; a link may have many different special alternating presentations.)

\subsection*{Log-concavity from Gr\"obner geometry}
We briefly summarize the way in which Gr\"obner geometry enters our story. First, the ideals $I_{L(\mathcal A)}$ defining $X_{\mathcal A}$ (and $\widehat X_G$) come from toric geometry \cite{sturmfels91,sturmfels96}. In \cite{lp13}, $f_{\mathcal A}(t)$ is built out of a \emph{fine mixed subdivision} of the zonotope spanned by $\mathcal A$. We apply the \emph{Cayley trick} \cite{hrs00} to turn this mixed subdivision into a triangulation of the \emph{Lawrence polytope} $L(\mathcal A)$ of $\mathcal A$. The dictionary between polyhedral and toric geometry produces, out of a polytope equipped with a triangulation, a toric ideal $I_{L(\mathcal A)}$ equipped with a Gr\"obner degeneration \cite{sturmfels96}.

More specifically, when $\mathcal A$ is rational and flat, the technology of toric geometry produces
\[
X_{\mathcal A}\subseteq \PP_x^{N-1}\times\PP_y^{N-1}
\]
equipped with a Gr\"obner degeneration
\[
X_{\mathcal A}\rightsquigarrow \init_\rho(X_{\mathcal A}), \quad \textup{ where } \quad (\init_\rho(X_{\mathcal A}))_\red = \bigcup_{B\in \mathcal B(\mathcal A)} H_B,
\]
in which the scheme-theoretic multiplicity of $H_B$ in the limit equals $c\cdot \Vol(\Pi_B)$ for some scalar $c$ independent of $B$, and the class of $H_B$ in $A^*(\PP_x^{N-1}\times \PP_y^{N-1})$ is
\[
[H_B] = q^{N-n-\ext_\rho(B)}t^{\ext_\rho(B)},
\]
where $q$ and $t$ are the hyperplane classes of $\PP_x^{N-1}$ and $\PP_y^{N-1}$ respectively. Consequently, the Chow class
\[
[X_{\mathcal A}] = \sum_B c\cdot \Vol(\Pi_B)[H_B]
\]
is a scalar multiple of the homogenization of $f_{\mathcal A}(t)$. These considerations imply the following result.
\begin{thm}
\label{thm:log-concave}
Let $\mathcal A$ be any full-dimensional flat vector arrangement. Then the coefficients of $f_{\mathcal A}(t)$ are log-concave.
\end{thm}
The geometric argument above requires $\mathcal A$ to be rational, but one can approximate an arbitrary arrangement $\mathcal A$ by rational arrangements.

Gao--Yuan have also independently proved Theorem~\ref{thm:log-concave} in \cite{gy26b}, building on their previous work \cite{gy26a}. We circulated a preliminary draft of this paper, which had a proof of Theorem~\ref{thm:log-concave}, about three weeks before \cite{gy26b} appeared on the arXiv. Around this time, we were just beginning to realize that multiplicity-freeness of $\mathfrak S_G$ bodes well for the study of $\mathfrak G_G$, which resulted in Theorems~\ref{thm:G_G-top} and~\ref{thm:G_G-extra} (along with an arXiv-upload delay).

\subsection*{Roadmap of the paper} We discuss background on external semi-activity, mixed subdivisions, degenerations of toric varieties, twisted $K$-polynomials, and M-convexity in Section~\ref{sec:background}. In Section~\ref{sec:cayley} we apply the Cayley trick to Li--Postnikov's zonotopal tiling and describe the corresponding triangulation of the Lawrence polytope. In Section~\ref{sec:graphical} we study the variety $\widehat X_G$ and the degeneration $\widehat X_G\rightsquigarrow \widehat X_0$ to prove Theorem~\ref{thm:S_G}. In Section~\ref{sec:G_G} we study the twisted $K$-polynomial $\mathfrak G_G(\mathbf z)$ and prove Theorems~\ref{thm:G_G-top} and~\ref{thm:G_G-extra}. To finish, in Section~\ref{sec:log-concave} we prove Theorem~\ref{thm:log-concave}.

\subsection*{Use of AI} 
I learnt the reference \cite{sturmfels96} from ChatGPT 5.6 Sol when I asked it to explain whether the Stanley--Reisner ring of a triangulated polytope is the Gr\"obner limit of something. I used Claude to check for errors. Everything else in this paper is human-made.

\subsection*{Acknowledgements}
I thank Karola M\'esz\'aros for introducing me to this subject, for many patient explanations and fruitful discussions about $f_{\mathcal A}(t)$ and the super cool polytope $Z_G^-$, and incidentally also for teaching me the Cayley trick nine years ago. I thank Alexander Postnikov for many insightful discussions about zonotopes and tilings, and for explaining results in \cite{lp13} and \cite{ap10,hr11,lenz14}. I thank Daoji Huang and Allen Knutson for inspiring discussions about Gr\"obner geometry. I thank Daoji Huang, Matt Larson, and Ada Stelzer for valuable feedback on a previous version of this paper and for telling me about a lot of useful context. Finally, I thank Cameron Chang, Jack Chen-An Chou, Alexandra Pevzner, Saskia Solotko, Ada Stelzer, and Alan Yan for helpful conversations.
\section{Background}
\label{sec:background}

\subsection*{Regular mixed subdivisions and the Cayley trick}
We follow the conventions in \cite{sturmfels96,hrs00,postnikov09}. In order to describe the Cayley trick, we recall some generalities on subdivisions.

\begin{defn}
A \demph{subdivision} of a polytope $P = \mathrm{conv}\{p_1, \dots, p_N\} \subset \RR^n$ is a decomposition of $P$ into a union of cells of the same dimension as $P$, such that each cell is the convex hull of some subset of vertices of $P$, and that any two cells intersect in their common face. Subdivisions are partially ordered by refinement, and a minimal element of this partial order is called a \demph{triangulation}.
\end{defn}

A good source of subdivisions of $P$ comes from the following construction. Let us say that a \demph{height vector} for $P = \mathrm{conv}\{p_1, \dots, p_N\}$ is a vector $h = (h_1, \dots, h_N) \in \RR^N$.

\begin{defn}
Let $P = \mathrm{conv}\{p_1, \dots, p_N\} \subset \RR^n$ be a polytope and let $h = (h_1, \dots, h_N) \in \RR^N$ be a height vector for $P$. Lift $p_i \in \RR^n$ to $\tilde p_i \colonequals (p_i, h_i) \in \RR^{n+1}$, and lift $P$ to $\tilde P \colonequals \mathrm{conv}\{\tilde p_i\}$. A \demph{lower face} $F$ is a face of $\tilde P$ that is ``visible from below'': there exists an affine-linear functional $f \colon \RR^n \to \RR$ so that $f(p_i) = h_i$ for $p_i \in F$ and $f(p_i) < h_i$ for $p_i\not\in F$. The \demph{regular subdivision} of $P$ induced by $h$ is the subdivision obtained by projecting the lower faces of $\tilde P$ back to $P$.
\end{defn}

There is also a notion of subdivision with respect to Minkowski sums.

\begin{defn}
A \demph{mixed subdivision} of a Minkowski sum of polytopes $P = P^{(1)} + \dots + P^{(m)}$ is a decomposition of $P$ into a union of \demph{Minkowski cells}, i.e., into a union of polytopes $C = C^{(1)} + \dots + C^{(m)}$ where each $C^{(i)}$ is a convex hull of some subset of vertices of $P^{(i)}$, so that any two cells intersect in their common face. Mixed subdivisions are partially ordered by refinement, and a minimal element of this partial order is called a \demph{fine mixed subdivision}.
\end{defn}

\begin{lem}[{\cite[Proposition 2.3]{santos05}, \cite[Lemma 14.2]{postnikov09}}]
\label{lem:FINE}
A mixed subdivision is fine if and only if for each mixed cell $C = C^{(1)} + \dots + C^{(m)}$ in this subdivision, all $C^{(i)}$ are simplices and $\dim(C^{(1)}) + \dots + \dim(C^{(m)}) = \dim(C)$, i.e., the simplices $C^{(i)}$ span independent affine subspaces.
\end{lem}

\begin{defn}
Let $P^{(1)}, \dots, P^{(m)}\subset\RR^n$ be polytopes, let $P = P^{(1)} + \dots + P^{(m)}$ denote their Minkowski sum, and pick a height vector $h^{(j)}$ for each $P^{(j)}$. As before, lift the vertices $p_i^{(j)}$ of $P^{(j)}$ to $\tilde p_i^{(j)} \colonequals (p_i^{(j)}, h_i^{(j)}) \in \RR^{n+1}$, and lift $P^{(j)}$ to $\tilde P^{(j)} \colonequals \mathrm{conv}\{\tilde p^{(j)}_i\}\subset \RR^{n+1}$. The \demph{regular mixed subdivision} of $P$ induced by the heights $h^{(j)}$ is the mixed subdivision obtained by projecting the lower faces of $\tilde P^{(1)} + \dots + \tilde P^{(m)}$ back onto $P$.
\end{defn}

Concretely, the cells $C = C^{(1)} + \dots + C^{(m)}$ of the regular mixed subdivision of $P = P^{(1)} + \dots + P^{(m)}$ induced by the $h^{(j)}$ are those for which there exists a linear functional $f\colon \RR^n\to\RR$, and constants $c_1,\dots, c_m$, so that, for each $j$,
\[
f(p_i^{(j)}) = h_i^{(j)} + c_j \quad \textup{ if } p_i^{(j)} \in C^{(j)}, \qquad\textup{ and } \qquad f(p_i^{(j)}) < h_i^{(j)} + c_j \quad \textup{ if } p_i^{(j)} \not\in C^{(j)}.
\]
This is analogous to the description of lower faces of ordinary subdivisions of polytopes in terms of affine-linear functionals.

The Cayley trick turns mixed subdivisions of $P = P_1 + \dots + P_N\subset\RR^n$ into ordinary subdivisions of a larger polytope
\[
\Cay(P_1, \dots, P_N) \colonequals \conv\left(\bigcup_{i=1}^N P_i \times \{\mathbf 1_i\}\right)\subset\RR^n\times\RR^N.
\]
called the \demph{Cayley polytope}. This idea was used \cite{santos05} to study triangulations of products of simplices.

The Cayley trick is captured by the following two results.

\begin{prop}[{\cite[Lemma 3.2]{hrs00}, \cite[Lemma 14.4]{postnikov09}}]
\label{prop:cayley-polytope}
For any $\lambda_1, \dots, \lambda_N \geq 0$ satisfying $\lambda_1 +\dots+\lambda_N = 1$, the intersection of $\Cay(P_1, \dots,P_N)$ with the affine subspace $\RR^n\times \{\lambda\}$ is exactly the scaled Minkowski sum $\lambda_1 P_1 + \dots + \lambda_N P_N$.
\end{prop}

(For a scalar $\lambda > 0$ and a polytope $P\subset\RR^n$, we write $\lambda P$ for the polytope obtained by dilating $P$ by a factor of $\lambda$.)

\begin{thm}[{\cite[Propositions 3.3, 3.4]{hrs00}, \cite[Lemma 14.5]{postnikov09}}]
\label{thm:cayley-trick}
For strictly positive $\lambda_1, \dots, \lambda_N > 0$ satisfying $\lambda_1 +\dots+\lambda_N = 1$, the intersection of a subdivision of $\Cay(P_1, \dots, P_N)$ with the affine subspace $\RR^n\times\{\lambda\}$ gives a mixed subdivision of $\lambda_1 P_1 + \dots + \lambda_N P_N$. This process gives a poset isomorphism between subdivisions of $\Cay(P_1, \dots, P_N)$ and mixed subdivisions of $\lambda_1P_1 + \dots + \lambda_NP_N$.
\end{thm}

\subsection*{External semi-activity and zonotopes}
Let $\mathcal A = \{a_1, \dots, a_N\}$ be a collection of vectors. The \demph{matroid} $M[\mathcal A]$ associated to a vector configuration $\mathcal A = \{a_1, \dots, a_N\}$ is defined on the ground set $[N]$, and its bases are subsets $\{i_1, \dots, i_n\} \subset[N]$ such that $\{a_{i_1}, \dots, a_{i_n}\}$ forms a basis for the span of the vectors in $\mathcal A$. We write $\mathcal B(\mathcal A)$ for the set of bases of $\mathcal A$.

For any basis $B$ and any $j\not\in B$, the set $B\cup\{j\}$ contains a unique circuit $C = C(B,j)$ called the fundamental circuit of $j$ with respect to $B$. Conversely, every circuit $C$ is $C(B,j)$ for some $j$ and some basis $B$ not containing $j$.

For any circuit $C$ of $M[\mathcal A]$, there is a unique, up to rescaling, nontrivial linear dependence
\begin{equation}
\label{eqn:dependence}
\sum_{i\in C}\lambda_i^{(C)} a_i = 0.
\end{equation}

Vectors satisfying~\eqref{eqn:dependence} are called \demph{circuit vectors}. When a circuit $C$ is viewed as the fundamental circuit $C(B,j)$, the unique vector $\lambda^{(C)}$ satisfying Equation~\eqref{eqn:dependence} with $\lambda^{(C)}_j = 1$ is called a \demph{normalized circuit vector} and shall be denoted $\lambda^{(C,j)}$. A framing vector $\rho \in \RR^N$ is called \demph{generic} for $\mathcal A$ if $\langle \rho, \lambda^{(C)}\rangle \neq 0$ for all circuits $C$ and if its coordinates are linearly independent over $\QQ$.

\begin{defn}
For $B \in \mathcal B(\mathcal A)$, an element $j \not\in B$ is \demph{externally semi-active} with respect to $\rho$ if, in the unique circuit $C\subset B\cup\{j\}$, the inequality $\langle \rho, \lambda^{(C,j)}\rangle > 0$ holds. We write $\Ext_\rho(B)$ for the set of externally semi-active elements and $\ext_\rho(B) \colonequals |\Ext_\rho(B)|$.
\end{defn}
\begin{defn}[{\cite{lp13}}]
The \demph{generalized Alexander polynomial} is
\[
f_{\mathcal A}(t) \colonequals \sum_{B\in\mathcal B(\mathcal A)}\Vol(\Pi_B)t^{\ext_\rho(B)},
\]
where
\[
\Pi_B \colonequals \sum_{i \in B}[0,a_i]
\]
is the parallelepiped spanned by $B$.
\end{defn}

We say that $\mathcal A$ is \demph{unimodular} if $\mathrm{span}_\ZZ\{a_i \colon i \in B\} = \mathrm{span}_\ZZ(\mathcal A)$.

It is known (\cite{lp13}, see also \cite[Theorem 3.4]{kmp25}) that $f_{\mathcal A}(t)$ is independent of choice of generic framing $\rho$, and the results of this paper give an independent proof of this fact. (The choice of $\rho$ affects the Gr\"obner degeneration of $X_{\mathcal A}$, but not the variety itself.)

The \demph{zonotope} $Z_{\mathcal A}$ spanned by $\mathcal A$ is the Minkowski sum
\[
Z_{\mathcal A} \colonequals \sum_{i=1}^N[0,a_i].
\]
The zonotope spanned by an arrangement of linearly independent vectors $B = \{b_1, \dots, b_d\}$ is just the parallelepiped $\Pi_B$ whose edges are parallel to the $b_i$. By Lemma~\ref{lem:FINE}, a fine mixed subdivision of $Z_{\mathcal A}$ amounts to a tiling by translates of parallelepipeds $\Pi_B$, where $B$ is a basis of the matroid $M[\mathcal A]$.

The \demph{trimmed zonotope} is the polytope
\[
Z_{\mathcal A}^- \colonequals \{x \in \RR^{n}\colon x + \mathbf 1_i \in Z_{\mathcal A} \textup{ for all } i \in [n]\}.
\]

\subsection*{Degenerations of toric varieties}
By work of Sturmfels \cite{sturmfels96,sturmfels91}, Gr\"obner degenerations of toric varieties are closely related to regular triangulations. Let $\mathcal A = \{a_1, \dots, a_N\}$ be a \demph{rational} arrangement, meaning that $a_i \in \QQ^n$ for all $i$, and assume that $\mathcal A$ is flat. Let 
\[
I_{\mathcal A} = \left\langle \mathbf x^\alpha - \mathbf x^\beta \colon \alpha,\beta\in\NN^N \textup{ and } \sum_{i=1}^N \alpha_i a_i = \sum_{i=1}^N\beta_i a_i\right\rangle \subseteq \CC[x_1, \dots, x_N]
\]
be the corresponding toric ideal.

Given a weight $\omega\in\RR^m$, define the \demph{initial form} of a polynomial $f = \sum c_\alpha \mathbf x^\alpha$ to be the sum of all terms $c_\alpha \mathbf x^\alpha$ such that the inner product $\langle \omega, \alpha\rangle$ is maximal. Define the \demph{initial ideal} of $I_{\mathcal A}$ with respect to $\omega$ by
\[
\init_\omega(I_{\mathcal A}) \colonequals \langle \init_\omega(f) \colon f \in I_{\mathcal A}\rangle.
\]
\begin{prop}[{\cite[Corollary 1.10]{sturmfels96}}]
\label{prop:omega}
If $\omega\geq 0$ and $\init_\omega(I_{\mathcal A})$ is a monomial ideal, then there exists a term order $\prec_\omega$ so that $\init_\omega(I_{\mathcal A}) = \init_{\prec_\omega}(I_{\mathcal A})$.
\end{prop}

The following results were originally stated in terms of initial ideals with respect to term orders, but using Proposition~\ref{prop:omega} we rephrase them in terms of initial ideals of the form $\init_\omega$.

\begin{thm}[{\cite[Corollary 8.4]{sturmfels96}}]
\label{thm:special-reduced}
If $\omega \geq 0$ and $\init_\omega(I_{\mathcal A})$ is a monomial ideal, then the regular subdivision $\Delta_\omega$ of $\mathrm{conv}(\mathcal A)$ is a triangulation and
\[
\sqrt{\init_\omega(I_{\mathcal A})} = \bigcap_{\substack{\sigma \in \Delta_\omega\\\sigma\textup{ maximal}}} H_\sigma, \qquad H_\sigma \colonequals \langle x_i \colon i \not \in \sigma\rangle.
\]
\end{thm}

Recall that the \demph{normalized volume} of a simplex $\sigma$ is defined to be the index $[\Lambda_{\mathcal A} : \Lambda_\sigma]$, where $\Lambda_{\mathcal A}$ is the $\ZZ$-lattice generated by vectors in $\mathcal A$ and $\Lambda_\sigma$ is the sublattice generated by vertices of $\sigma$.

\begin{thm}[{\cite[Theorem 8.8]{sturmfels96}}]
\label{thm:special-multiplicity}
In the situation of Theorem~\ref{thm:special-reduced}, the multiplicity of $H_\sigma$ in $\init_\omega(I_{\mathcal A})$ is the normalized volume $\Vol_{\mathrm{norm}}(\sigma)$ of the simplex $\sigma$.
\end{thm}

\begin{thm}[{\cite[Corollary 8.9]{sturmfels96}}]
\label{thm:special-radical}
In the situation of Theorem~\ref{thm:special-reduced}, initial ideal $\init_\omega(I_{\mathcal A})$ is square-free if and only if $\Vol_{\mathrm{norm}}(\sigma) = 1$ for all maximal simplices $\sigma$. In this case, $\sqrt{\init_\omega(I_{\mathcal A})} = \init_\omega(I_{\mathcal A})$.
\end{thm}

\subsection*{Twisted $K$-polynomials and Stanley--Reisner rings}
Let $R = \CC[x_1, \dots, x_n]$ be a $\ZZ^k$-graded ring, and set $\mathbf d_i \colonequals \deg(x_i) \in \ZZ^k$. Assume (as will be the case for everything in this paper) that $R$ is positively graded in the sense that $\mathbf d_i \in \ZZ_{\geq 0}^k \setminus\{\mathbf 0\}$. Let $I$ be a homogeneous ideal, so that $R/I$ is again $\ZZ^k$-graded. The \demph{Hilbert series} of $R/I$ is the generating function
\[
\Hilb(R/I;\mathbf t) \colonequals \sum_{\mathbf a \in \ZZ^k} \dim_\CC((R/I)_{\mathbf a})\mathbf t^{\mathbf a}
\]
and the \demph{$K$-polynomial} is
\[
\mathcal K(R/I; \mathbf t) \colonequals \Hilb(R/I;\mathbf t) \cdot \prod_{i=1}^n (1 - \mathbf t^{\mathbf d_i}).
\]
(The quantity $\mathcal K(R/I;\mathbf t)$ is always a polynomial, and records terms in the minimal free resolution of $R/I$.)

\begin{defn}
The \demph{twisted $K$-polynomial} of $R/I$ is the polynomial $\mathfrak G(R/I; \mathbf t) \colonequals \mathcal K(R/I; \mathbf 1- \mathbf t)$, i.e., the polynomial obtained from $\mathcal K(R/I;\mathbf t)$ by applying the substitution $t_i \mapsto 1 - t_i$.
\end{defn}

We remark (c.f.\ \cite[Theorem 8.36]{ms05}) that Gr\"obner degeneration preserves $\Hilb(\cdot)$ and hence preserves $\mathcal K(\cdot)$ and $\mathfrak G(\cdot)$.

Let $\Delta$ be a simplicial complex on vertex set $[n]$. The \demph{Stanley--Reisner ideal} $I_\Delta$ of is the ideal 
\[
I_\Delta \overset{\rm def}= \left\langle \prod_{i \in F} x_i \colon F\textup{ is not a face of } \Delta\right\rangle \subset \CC[x_1, \dots, x_n].
\]
\begin{rem}
\label{rem:init-is-SR}
If $\omega \geq 0$ and $\init_\omega(I_{\mathcal A})$ is a monomial ideal, the radical $\sqrt{\init_\omega(I_{\mathcal A})}$ is the Stanley--Reisner ideal of the triangulation $\Delta_\omega$ of $\mathrm{conv}(\mathcal A)$, viewed as a simplicial complex. (This follows immediately from Theorem~\ref{thm:special-reduced}.)
\end{rem}
The $\mathcal K$-polynomials of Stanley--Reisner rings $R/I_\Delta$ can be computed by a theorem of Knutson--Miller--Yong when $\Delta$ is homeomorphic to a ball. More precisely, if $\Delta$ is a ball, let us say that a face $F$ of $\Delta$ is a \demph{boundary face} if $F$ is contained in a face of codimension $1$ which lies in exactly one maximal face of $\Delta$.

\begin{thm}[{\cite[Proposition 4.2]{kmy08}}]
\label{thm:kmy}
Let $\Delta$ be a simplicial complex homeomorphic to a ball. Then
\[
\mathcal K(R/I_\Delta; \mathbf t) = \sum_F (-1)^{\codim F} \prod_{i\in [n]\setminus F}(1 - \mathbf t^{\mathbf d_i}),
\]
where the sum runs over the interior faces $F$ of $\Delta$.
\end{thm}

\subsection*{M-convexity and generalized polymatroids}

Let $f(x_1, \dots, x_n)$ be a polynomial. The \demph{support} of $f$ is set $\{\alpha \in \ZZ_{\geq 0}^n \colon \mathbf x^\alpha \textup{ appears in } f\}$ of exponent vectors of monomials appearing with nonzero coefficient in $f$. The \demph{Newton polytope} of $f$, denoted $\Newton(f)$, is the convex hull $\conv(\supp(f))$ of its support. 

Monical--Tokcan--Yong \cite{mty19} initiated the systematic study of polynomials $f$ with the \demph{saturated Newton polytope} property, abbreviated SNP: we say $f$ satisfies SNP if $\supp(f) = \Newton(f)\cap\ZZ^n$.

A set $S\subset\ZZ^n$ is \demph{M-convex} if for any $\alpha,\beta\in S$ and any $i\in[n]$ such that $\alpha_i > \beta_i$, there exists $j\in[n]$ such that $\alpha_j < \beta_j$ and $\alpha - \mathbf e_i + \mathbf e_j \in S$ and $\beta + \mathbf e_i - \mathbf e_j \in S$.

A function $z\colon 2^{[n]}\to\RR$ is called \demph{submodular} if
\[
z(I) + z(J) \geq z(I\cup J) + z(I\cap J)\qquad\textup{ for all } I,J\subseteq[n].
\]
A \demph{generalized permutahedron} $P\subset\RR^n$ is a polytope defined by
\[
P = \left\{x\in\RR^n\colon \sum_{i\in I}x_i \leq z(I) \textup{ for all } I\subseteq[n]\textup{ and } \sum_{i=1}^n x_i = z([n])\right\}
\]
for some submodular function $z\colon 2^{[n]}\to\RR$ satisfying $z(\emptyset) = 0$. A set is M-convex if and only if it is the set of integer points of an integral generalized permutahedron \cite[Theorem 4.15]{murota03}. We refer to \cite{schrijver03} and \cite{postnikov09} for background on generalized permutahedra.

We need the following strengthened form of the integer decomposition property of generalized permutahedra.

\begin{lem}[{\cite[Theorem 44.6, Corollary 46.2c]{schrijver03}}]
\label{lem:gp-minkowski}
The Minkowski sum $P_1 + \dots + P_k$ of generalized permutahedra is a generalized permutahedron. Furthermore, if each $P_i$ is integral,
\[
\left(\sum_{j=1}^k P_j\right) \cap \ZZ^n = \sum_{j=1}^k(P_j\cap\ZZ^n),
\]
that is to say, any integer point in a Minkowski sum of generalized permutahedra can be written as a sum of integer points of summands.
\end{lem}

\subsection*{Multiplicity-freeness}
Our next goal is to recall \cite[Theorem A]{ccmm22} and \cite[Theorem 3]{knutson09}.

Let $S\subset \ZZ_{\geq 0}^n$ be a finite set, and set $N\colonequals \max_{s\in S} |s|$ (where $|s|$ denotes the sum of the coordinates of $s \in S$). The \demph{homogenization} $\widetilde S$ of $S$ is the set
\[
\widetilde S \colonequals \{(s,N-|s|)\colon s \in S\}\subset \ZZ_{\geq 0}^{n+1}\}
\]

Following \cite[Defn 5.5]{ccmm22}, a \demph{generalized polymatroid} is a set $S\subset \ZZ_{\geq 0}^n$ so that its homogenization $\widetilde S$ is M-convex.

We say that a subvariety $X \subset \PP \colonequals \PP^{m_1} \times \dots \times \PP^{m_p}$ is \demph{multiplicity-free} if the corresponding class $[X] \in A^*(\PP) \cong \ZZ[z_1, \dots, z_p]/\langle z_1^{m_1+1}, \dots, z_p^{m_p+1}\rangle$ is a multiplicity-free polynomial, i.e., has coefficients zero or one.

\begin{thm}[{\cite[Theorem A]{ccmm22}}]
\label{thm:generalized-polymatroid-support}
Let $X \subset \PP$ be a multiplicity-free variety. Then, the support of the twisted $K$-polynomial $\mathfrak G(X;\mathbf z)$ is a generalized polymatroid.
\end{thm}
Theorem~\ref{thm:generalized-polymatroid-support} implies, in particular, that the support of $\mathfrak G(X;\mathbf z)$ (and its homogenization) is saturated.

Castillo--Cid-Ruiz--Mohammadi--Monta\~no also give an explicit description of the support of $\mathfrak G(X;\mathbf z)$ in \cite[Theorem A]{ccmm22}, though we will not need it in this paper.

In order to simplify exposition, we will state \cite[Theorem 3]{knutson09} specialized to multiplicity-free subvarieties of products of projective spaces. (We follow the exposition in \cite{knutson09}.)

The Schubert subvarieties of $\PP$ are products $X_{\mathbf p} \colonequals \prod_i \PP^{m_i - p_i}$ of projective linear subspaces, and their Chow class is $[X_{\mathbf p}] = \mathbf z^{\mathbf p}$. A Schubert variety $X_{\mathbf p}$ contains another Schubert variety $X_{\mathbf q}$ if $\mathbf p \leq \mathbf q$ (componentwise comparison).

Let $X\subset \PP$ be a multiplicity free variety, and write
\[
[X] = \sum_{\mathbf p \in S} \mathbf z^{\mathbf p} \in \ZZ[z_1, \dots, z_p]/\langle z_1^{m_1+1}, \dots, z_p^{m_p+1}\rangle, \qquad S \subseteq \{\mathbf p \colon 0 \leq p_i \leq m_i \textup{ for all } i\},
\]
so that $S$ is the support of the polynomial representing $[X]$. A theorem of Brion \cite{brion03} guarantees that $X$ flatly degenerates to the reduced union $\bigcup_{\mathbf p \in S} X_{\mathbf p}$ of Schubert subvarieties of $\PP$. Write
\[
P_S \colonequals \{\mathbf q \in \ZZ_{\geq 0}^p \colon \mathbf q \leq \mathbf m \textup{ and } \mathbf q \geq \mathbf p \textup{ for some } \mathbf p \in S\};
\]
this poset is isomorphic to the poset of Schubert varieties $X_{\mathbf q} \subset \bigcup_{\mathbf p \in S} X_{\mathbf p}$ ordered by reverse containment.

Let us call the \demph{M\"obius function} of a poset $P$ to be the unique function $\mu\colon P\to\ZZ$ satisfying
\[
\sum_{\mathbf q\leq_P \mathbf p}\mu(\mathbf q) = 1, \qquad\textup{ for all } \mathbf p \in P.
\]

Then \cite[Theorem 3]{knutson09} can be stated as follows.
\begin{thm}
\label{thm:mobius}
Let $X$ be a multiplicity-free subvariety of $\PP$. Then
\[
\mathcal K(X;\mathbf z) = \sum_{\mathbf q \in P_S} \mu(\mathbf q) \mathcal K(X_{\mathbf q} ;\mathbf z),
\]
where $\mu$ is the M\"obius function of $P_S$.
\end{thm}
\begin{cor}
\label{cor:mobius-explicit}
Let $X$ be a multiplicity-free subvariety of $\PP$. Then
\[
\mathfrak G(X;\mathbf z) = \sum_{\mathbf q \in P_S} \mu(\mathbf q) \mathbf z^{\mathbf q}
\]
where $\mu$ is the M\"obius function of $P_S$.
\end{cor}
\begin{proof}
Substitute $z_i \mapsto 1 - z_i$ in Theorem~\ref{thm:mobius} and compute that
\[
\mathfrak G(X_{\mathbf q}; \mathbf z) = \mathbf z^{\mathbf q}.
\]
(The variety $X_{\mathbf q}$ is is a product of linear subspaces and, e.g., its Hilbert series can be computed explicitly.)
\end{proof}

\section{The Cayley trick}
\label{sec:cayley}
We recall a fine mixed subdivision of $Z_{\mathcal A}$ by translates of $\Pi_B$ (Theorem~\ref{thm:tiling}) studied by Li--Postnikov \cite{lp13}. The goal of this section is to apply the Cayley trick (Theorem~\ref{thm:cayley-trick}) to the tiling in Theorem~\ref{thm:tiling}. The resulting triangulated polytope gives a Gr\"obner degeneration of the toric ideal $I_{L(\mathcal A)}$ to a monomial ideal (Lemma~\ref{lem:init-is-monomial}).

\begin{thm}[{\cite{lp13}; see also \cite[Section 3]{kmp25}}]
\label{thm:tiling}
The vector $\rho$ determines a regular fine mixed subdivision by the tiles
\[
\mathcal E_{B,\rho}\colonequals \sum_{i\in B}[0,a_i] + \sum_{j\in\Ext_\rho(B)}a_j
\]
indexed by bases $B$ of $\mathcal A$.
\end{thm}
This result appears to be well-known in the literature, but we include a proof for completeness.
\begin{proof}
For each $i$, define the height function $h^{(i)}$ on the vertices of the $i$-th segment $[0,a_i]$ by
\begin{align*}
h^{(i)}(0)&\colonequals 0,\\
h^{(i)}(a_i)&\colonequals - \rho_i.
\end{align*}
We claim that the regular mixed subdivision defined by the $h^{(i)}$ is the desired tiling.

The lower face of the lift of $[0,a_i]$ with respect to a linear functional $\ell\colon \RR^n\to\RR$ is
\[
C_i(\ell) = \begin{cases} \{0\} &\textup{ if } \ell(a_i) < -\rho_i,\\ [0,a_i] &\textup{ if } \ell(a_i) = -\rho_i,\\ \{a_i\}&\textup{ if } \ell(a_i) > -\rho_i.
\end{cases}
\]
Given a linear functional $\ell$, let $S_\ell \colonequals \{i\in[N] \colon \ell(a_i) = -\rho_i\}$. If $S_\ell$ contained a circuit $C$, say with circuit vector $\lambda^{(C)}$, then the equalities $\ell(a_i) = -\rho_i$ for all $i \in C$ imply that
\begin{align*}
\langle \lambda^{(C)},\rho\rangle &= \sum_{i\in C}\lambda_i^{(C)}\rho_i
\\&=\ell\left(-\sum_{i\in C}\lambda_i^{(C)}a_i\right)
\\&= \ell(0)
\\&= 0,
\end{align*}
contradictory to assumption on $\rho$. Thus, $S_\ell$ must be an independent set, and the corresponding tile has dimension $|S_\ell|$. In particular, the regular mixed subdivision induced by $h^{(i)}$ is fine.

It remains to characterize the maximal cells of this subdivision as the tiles $\mathcal E_{B,\rho}$. To this end, first observe that given a basis $B$ there exists a linear functional $f$ satisfying $f(a_i) = \rho_i$ for all $i \in B$, i.e., $B = S_{-f}$ for some $f$. Then, for any $j \not \in B$, let $C$ denote the circuit in $B\cup\{j\}$ and recall that
\[
\begin{cases} \langle \lambda^{(C,j)},\rho\rangle > 0 &\textup{ if } j\in \Ext_\rho(B),\\
\langle \lambda^{(C,j)},\rho\rangle < 0 &\textup{ if } j\not\in \Ext_\rho(B).
\end{cases}
\]
From the fact that
\begin{align*}
0 = \sum_{i\in C}\lambda_i^{(C,j)} f(a_i) = \sum_{i \in C\setminus \{j\}} \lambda_i^{(C,j)} \rho_i + \lambda_j^{(C,j)} f(a_j),
\end{align*}
we deduce that
\begin{equation}
\label{eqn:again}
\begin{cases} f(a_j) < \rho_j &\textup{ if } j \in \Ext_\rho(B),\\ f(a_j) > \rho_j &\textup{ if } j\not\in \Ext_\rho(B).\end{cases}
\end{equation}
Hence, the lower face of the lift of $[0,a_i]$ for the functional $\ell \colonequals -f$ is
\[
C_i(\ell) = \begin{cases} [0,a_i] &\textup{ if } i \in B,\\ \{a_i\} &\textup{ if } i \in \Ext_\rho(B),\\ \{0\} &\textup{ if } i\not\in B\cup \Ext_\rho(B),\end{cases}
\]
and so the resulting Minkowski cell in the mixed subdivision is
\[
\sum_{i=1}^N C_i(\ell) = \mathcal E_{B,\rho}.\qedhere
\]
\end{proof}

\begin{defn}[{cf.\ \cite[\S 4.2]{hrs00}}]
Let $\mathcal A = \{a_1, \dots, a_N\}\subset\RR^n$ denote a vector arrangement, and define vectors $b_i\colonequals (a_i, \mathbf 1_i) \in \RR^{n+N}$ and $c_i\colonequals (0,\mathbf 1_i) \in \RR^{n+N}$. The \demph{Lawrence configuration} is
\[
L(\mathcal A) = \{b_i, c_i \colon i\in[N]\},
\]
and the \demph{Lawrence polytope} is $\conv(L(\mathcal A))$.
\end{defn}
By construction, $L(\mathcal A)$ is the Cayley polytope of the zonotope $Z_{\mathcal A}$. Every Lawrence configuration is flat, with respect to the linear functional $(a,z)\mapsto|z|$.

\begin{prop}
\label{prop:triangulation}
Under the Cayley trick, the tiling of Theorem~\ref{thm:tiling} is sent to a regular triangulation $\mathcal T_\rho$ of $\conv(L(\mathcal A))$, induced by the height function $h(b_i) = -\rho_i$ and $h(c_i) = 0$. The simplices in $\mathcal T_\rho$ are
\[
T_{B,\rho} \colonequals \conv(\{b_i, c_i \colon i \in B\} \cup \{b_j \colon j \in \Ext_\rho(B)\} \cup \{c_j \colon j \not\in B\cup \Ext_\rho(B)\}).
\]
\end{prop}
\begin{proof}
We claim that it suffices to exhibit a linear functional whose lower face is $T_{B,\rho}$ for any basis $B$: The cell $T_{B,\rho}$ is itself a Cayley polytope and Proposition~\ref{prop:cayley-polytope} implies that 
\[
T_{B,\rho} \cap \left(\RR^n\times\left\{\left(\frac1N,\dots,\frac1N\right)\right\}\right) = \frac1N\mathcal E_{B,\rho}.
\]
Then Theorem~\ref{thm:cayley-trick} implies that $\mathcal T_\rho$ is a triangulation, as its restriction to $Z_{\mathcal A}$ is a fine mixed subdivision. The cells $T_{B,\rho}$ then exhaust $\conv(L(\mathcal A))$, as any other maximal cell $C$ is, again, a Cayley polytope and Proposition~\ref{prop:cayley-polytope} produces a tile $C\cap\RR^n\times\left\{\left(\frac1N,\dots,\frac1N\right)\right\}$ in the induced zonotopal tiling of $Z_{\mathcal A}$, contradicting the fact that the $\mathcal E_{B,\rho}$ exhaust $Z_{\mathcal A}$.

Given a basis $B$, let $f\colon \RR^n\to\RR$ denote the unique linear functional sending $f(a_i) = \rho_i$ for all $i \in B$. Recall from Equation~\eqref{eqn:again} that
\[
\begin{cases} f(a_j) < \rho_j &\textup{ if } j \in \Ext_\rho(B),\\ f(a_j) > \rho_j &\textup{ if } j\not\in \Ext_\rho(B).\end{cases}
\]
Define $\ell\colon \RR^n\times\RR^N \to \RR$ by the formula
\[
\ell(a,z) = -f(a) + \sum_{j\in\Ext_\rho(B)}(f(a_j)-\rho_j)z_j.
\]
Then observe that:
\begin{itemize}
\item For $i \in B$, we have
\[
\ell(a_i, \mathbf 1_i) = -\rho_i \qquad\textup{ and } \qquad \ell(0,\mathbf 1_i) = 0,
\]
\item For $j\in\Ext_\rho(B)$, we have 
\[
\ell(a_j,\mathbf 1_j) = -\rho_j \qquad\textup{ and } \qquad \ell(0,\mathbf 1_j) = f(a_j) - \rho_j < 0,
\]
\item For $j\not\in B\cup\Ext_\rho(B)$, we have
\[
\ell(a_j, \mathbf 1_j) = -f(a_j) < -\rho_j \qquad\textup{ and } \qquad \ell(0,\mathbf 1_j) = 0.
\]
\end{itemize}
Thus $\ell(b_i) \leq -\rho_i$ and $\ell(c_i)\leq 0$ for all $i$, and equality holds if and only if $b_i$ or $c_i$ is a vertex of $T_{B,\rho}$.
\end{proof}
\begin{rem}
\label{rem:minimal-simplicies}
The proof of Proposition~\ref{prop:triangulation} also implies the folklore observation that the poset of interior faces in the triangulation $\mathcal T_\rho$ of $\mathrm{conv}(L(\mathcal A))$ is isomorphic to the poset of interior faces in the fine mixed subdivision of Theorem~\ref{thm:tiling}. (This observation is also implicit in the proof of \cite[Lemma 3.2]{hrs00}.)
\end{rem}
For later use (cf.\ Theorem~\ref{thm:special-reduced}), we record the following corollary.
\begin{cor}
\label{cor:nonvertices}
The nonvertices of the simplex $T_{B,\rho}$ are $\{b_j \colon j \not\in B\cup\Ext_\rho(B)\} \cup \{c_j \colon j \in \Ext_\rho(B)\}$.
\end{cor}
\begin{proof}
Take the complement of the set of vertices of $T_{B,\rho}$, as described by Proposition~\ref{prop:triangulation}.
\end{proof}

We also record here a lemma for use in Section~\ref{sec:G_G}.

\begin{lem}
\label{lem:interior-necessary}
Suppose that $\mathrm{conv}(S)$ is an interior simplex in the triangulation $\mathcal T_\rho$. Then $|S \cap \{b_i, c_i\}| \geq 1$ for all $i$. In particular, the interior faces of $\mathcal T_\rho$ have dimension at least $N-1$.
\end{lem}
\begin{proof}
If $S\cap\{b_i, c_i\} = \emptyset$, then $\mathrm{conv}(S)$ lives in the supporting hyperplane $\{(a,z) \in \RR^n\times\RR^N\colon z_i = 0\}$ of $\mathrm{conv}(L(\mathcal A))$.

It follows that an interior face has at least $N$ vertices, and so $\dim(\mathrm{conv}(S)) \geq N-1$.
\end{proof}

Our next goal is to pass the triangulation $\mathcal T_\rho$ of $\mathrm{conv}(L(\mathcal A))$ through the toric dictionary. We work in the polynomial ring $R\colonequals \CC[x_1, \dots, x_N, y_1, \dots, y_N]$, corresponding to elements in $L(\mathcal A)$, where $\mathcal A$ is a rational arrangement. Recall that $I_{L(\mathcal A)}$ is the ideal generated by binomials of the form
\[
\mathbf x^\alpha \mathbf y^{\alpha'} - \mathbf x^\beta \mathbf y^{\beta'}, \qquad \textup{ where }\sum_{i=1}^N \alpha_i b_i + \sum_{i=1}^N \alpha_i' c_i = \sum_{i=1}^N \beta_i b_i + \sum_{i=1}^N \beta_i'c_i.
\]
Define the weight vector $L(\rho) \in \RR^{2N}$ by first choosing $M \in \RR_{>0}$ satisfying $M > \rho_i$ for all $i$, and then taking
\[
L(\rho)_i = \begin{cases} M - \rho_i &\textup{ if } i \in \{1,\dots,N\}, \\ M &\textup{ for } i \in \{N+1,\dots, 2N\},\end{cases}
\]
i.e., by giving weight $M - \rho_i$ to variables $x_i$ and weight $M$ to variables $y_i$. The regular subdivision of $\conv(L(\mathcal A))$ defined by the height function $L(\rho)$ is $\mathcal T_\rho$.
\begin{lem}
\label{lem:init-is-monomial}
The ideal $\init_{L(\rho)}(I_{L(\mathcal A)})$ is monomial.
\end{lem}
\begin{proof}
Let $<$ be a term order refining $L(\rho)$. The ideal $I_{L(\mathcal A)}$ is binomial, so the reduced Gr\"obner basis $G$ for $<$ consists of binomials.

Take any binomial in $I_{L(\mathcal A)}$, say with monomials $\mathbf x^\alpha\mathbf y^{\alpha'}$ and $\mathbf x^\beta\mathbf y^{\beta'}$. Then $|\alpha| + |\alpha'| = |\beta| + |\beta'|$ because $I_{L(\mathcal A)}$ is homogeneous (with respect to the usual degree $\deg(x_i) = \deg(y_i) = 1$). Thus, the $L(\rho)$-weights of $\mathbf x^\alpha\mathbf y^{\alpha'}$ and $\mathbf x^\beta\mathbf y^{\beta'}$ are equal to
\[
(|\alpha| + |\alpha'|)M - \langle \rho,\alpha\rangle \qquad\textup{ and } \qquad (|\beta| + |\beta'|)M - \langle \rho,\beta\rangle,
\]
and, since $\rho$ is generic, the fact that $\langle \rho,\alpha-\beta\rangle\neq 0$ implies that $\mathbf x^\alpha\mathbf y^{\alpha'}$ and $\mathbf x^\beta\mathbf y^{\beta'}$ have distinct $L(\rho)$-weight.

It follows that $\init_{L(\rho)}(g)=\init_<(g)$ for each $g\in G$, so $\init_<(I_{L(\mathcal A)})\subseteq\init_{L(\rho)}(I_{L(\mathcal A)})$. Both ideals have the same Hilbert function as $I_{L(\mathcal A)}$, so this inclusion is an equality, and in particular $\init_{L(\rho)}(I_{L(\mathcal A)})$ is monomial.
\end{proof}

Let $\Lambda_{\mathcal A} \colonequals \sum_{i=1}^N\ZZ a_i$ denote the $\ZZ$-lattice generated by $\mathcal A$.
\begin{lem}
\label{lem:volumes}
For any basis $B$, the equality
\[
\Vol(\Pi_B) = \Vol_{\mathrm{norm}}(T_{B,\rho}) \cdot \covol(\Lambda_{\mathcal A})
\]
holds.
\end{lem}
\begin{proof}
Observe that the lattice $\Lambda_{L(\mathcal A)}$ generated by $L(\mathcal A)$ equals $\Lambda_{\mathcal A} \oplus \ZZ^N$. Recall that normalized volume is defined to be $\Vol_{\mathrm{norm}}(T_{B,\rho}) = [\Lambda_{L(\mathcal A)} : \Lambda_{T_{B,\rho}}]$. The lattice $\Lambda_{T_{B,\rho}}$ is generated by 
\[
\{(a_i, 0) \colon i \in B\}\cup\{(0, \mathbf 1_i) \colon i \in B\}\cup\{(a_i, \mathbf 1_i) \colon i \in \Ext_\rho(B)\}\cup\{(0,\mathbf 1_i) \colon i \not\in B\cup\Ext_\rho(B)\}.
\]
For $i \in [N]$, let us write $v_i$ for the unique vector in the list above that has second coordinate equal to $\mathbf 1_i$. Then any vector $v = (a_1, \dots, a_n, z_1, \dots, z_N) \in \Lambda_{L(\mathcal A)}$ can be written as
\[
v = \sum_{i=1}^N z_iv_i + (a_1', \dots, a_n', 0,\dots,0), \qquad a' \in \Lambda_{\mathcal A}, \textup{ and } z_i \in \ZZ.
\]
There is a surjective map of lattices
\begin{align*}
\varphi \colon \Lambda_{\mathcal A} \oplus \ZZ^N &\twoheadrightarrow \Lambda_{\mathcal A}\\
(a,z)&\mapsto a'.
\end{align*}
An element $(a,z) \in \Lambda_{L(\mathcal A)}$ is sent to $\Lambda_B$ if and only if $(a,z) \in \Lambda_{T_{B,\rho}}$. Thus $\varphi$ descends to an isomorphism
\[
\frac{\Lambda_{\mathcal A}\oplus\ZZ^N}{\Lambda_{T_{B,\rho}}} \cong \frac{\Lambda_{\mathcal A}}{\Lambda_B}.
\]
The claim follows from the equalities
\begin{align*}
\Vol(\Pi_B) &= [\Lambda_{\mathcal A} : \Lambda_B] \cdot \covol(\Lambda_{\mathcal A}),\\
\Vol_{\mathrm{norm}}(T_{B,\rho}) &= [\Lambda_{\mathcal A} \oplus\ZZ^N : \Lambda_{T_{B,\rho}}].\qedhere
\end{align*}
\end{proof}

\section{Graphical arrangements and multidegree of $\widehat X_G$}
\label{sec:graphical}
In this section, we study the graphical arrangement $\mathcal A_G$. We find a $\ZZ^{V(G)}$-grading with respect to which the ideal $I_{L(\mathcal A_G)}$ is homogeneous and compute the Chow class of the corresponding multiprojective variety $\widehat{X}_{\mathcal A_G}$ in Theorem~\ref{thm:multidegree} and to relate it to the reciprocal of the integer point transform of the trimmed zonotope of $G$. This gives Theorem~\ref{thm:S_G}.

For an edge $e \in E(G)$, write $s(e)$ and $t(e)$ for the source and target vertices of $e$ respectively. Let $\mathcal A_G = \{\mathbf 1_{s(e)} - \mathbf 1_{t(e)} \colon e \in E(G)\} \subset\RR^{V(G)}$ denote the \demph{graphical arrangement} associated to $G$. A subset $S\subseteq E(G)$ of edges defines a basis $\{\mathbf 1_{s(e)} - \mathbf 1_{t(e)} \colon e \in S\}$ of $\mathcal A_G$ precisely when $S$ forms the set of edges of a tree $T$, and in this case $\Vol(\Pi_T) = 1$.

Let us emphasize that the Lawrence configuration $L(\mathcal A_G)$ is the vector configuration in $\RR^{V(G)} \times \RR^{E(G)}$ consisting of vectors $\{b_e, c_e \colon e \in E(G)\}$ where $b_e = \mathbf 1_{s(e)} - \mathbf 1_{t(e)} + \mathbf 1_e$ and $c_e = \mathbf 1_e$. The toric ideal $I_{L(\mathcal A_G)}$ sits in a polynomial ring $R = \CC[x_e, y_e]$ in which the variable $x_e$ corresponds to the vector $b_e$ and the variable $y_e$ corresponds to the vector $c_e$; this ideal is generated by equations
\begin{equation}
\label{eqn:ilag-gen}
\mathbf x^\alpha \mathbf y^{\alpha'} - \mathbf x^\beta\mathbf y^{\beta'}, \qquad\textup{ where } \sum_{e\in E(G)}\alpha_e b_e + \sum_{e\in E(G)}\alpha_e' c_e = \sum_{e\in E(G)}\beta_e b_e + \sum_{e\in E(G)}\beta_e' c_e.
\end{equation}
\begin{lem}
\label{lem:multihomog}
The ideal $I_{L(\mathcal A_G)}$ is homogeneous with respect to the $\ZZ^{V(G)}$-valued grading defined by $\mdeg(x_e) \colonequals \mathbf 1_{s(e)}$ and $\mdeg(y_e) \colonequals \mathbf 1_{t(e)}$.
\end{lem}

\begin{proof}
Restricting the condition on $(\alpha,\alpha',\beta,\beta')$ in~\eqref{eqn:ilag-gen} to the coordinates in $\RR^{V(G)}$ gives
\begin{equation}
\label{eqn:multihomog-1}
\sum_{e\in E(G)}\alpha_e(\mathbf 1_{s(e)} - \mathbf 1_{t(e)}) = \sum_{e\in E(G)}\beta_e(\mathbf 1_{s(e)} - \mathbf 1_{t(e)}),
\end{equation}
while restricting to the coordinate indexed by $e$ in $\RR^{E(G)}$ gives $\alpha_e + \alpha_e' = \beta_e + \beta_e'$ for all $e \in E(G)$. Hence,
\begin{equation}
\label{eqn:multihomog-2}
\sum_{e\in E(G)}(\alpha_e + \alpha_e')\mathbf 1_{t(e)} = \sum_{e\in E(G)}(\beta_e + \beta_e')\mathbf 1_{t(e)},
\end{equation}
and adding Equations~\eqref{eqn:multihomog-1} and~\eqref{eqn:multihomog-2} yields
\[
\sum_{e\in E(G)}\alpha_e\mathbf 1_{s(e)} + \sum_{e\in E(G)}\alpha_e'\mathbf 1_{t(e)} = \sum_{e\in E(G)}\beta_e\mathbf 1_{s(e)} + \sum_{e\in E(G)}\beta_e'\mathbf 1_{t(e)}.
\]
Thus the difference $\mathbf x^\alpha \mathbf y^{\alpha'} - \mathbf x^\beta\mathbf y^{\beta'}$ is homogeneous with respect to $\mdeg$.
\end{proof}
\begin{defn}
Let
\[
\widehat X_G \subseteq \prod_{v \in V(G)} \PP_v^{\deg(v) - 1},
\]
denote the multiprojective variety associated to the homogeneous ideal $I_{L(\mathcal A)}$. 
\end{defn}
We emphasize that $\widehat X_G$ lives in a projectivization of a $2|E(G)|$-dimensional affine space with variables $\{x_e, y_e\colon e \in E(G)\}$. The homogeneous coordinates of $\PP_v^{\deg(v)-1}$ are the $\deg(v)$ many variables whose multiweight $\mdeg$ equals $\mathbf 1_v \in \ZZ^{V(G)}$, namely, $x_e$ for $s(e) = v$ and $y_e$ for $t(e) = v$.

\begin{thm}
\label{thm:multidegree}
For $v \in V(G)$, let $z_v$ denote the hyperplane class of $\PP_v^{\deg(v) - 1}$. Then
\[
[\widehat X_G] = \sum_T\mathbf z(T), \qquad \mathbf z(T) \colonequals \left(\prod_{e\not\in T\cup\Ext_\rho(T)}z_{s(e)}\right)\left(\prod_{e\in\Ext_\rho(T)}z_{t(e)}\right),
\]
where the sum runs over spanning trees $T$ of $G$.
\end{thm}

\begin{proof}
By Lemma~\ref{lem:init-is-monomial}, $\init_{L(\rho)}(I_{L(\mathcal A_G)})$ is monomial, and Proposition~\ref{prop:omega} gives a term order $\prec_\omega$ so that $\init_{L(\rho)}(I_{L(\mathcal A_G)}) = \init_{\prec_\omega}(I_{L(\mathcal A_G)})$. Hence, there is a Gr\"obner degeneration
\[
\widehat X_G \rightsquigarrow \widehat X_0
\]
where $\widehat X_0$ is the multiprojective scheme associated to $\init_{L(\rho)}(I_{L(\mathcal A_G)})$. Since $\mathcal A_G$ is unimodular, Lemma~\ref{lem:volumes} implies that all simplices in $\mathcal T_\rho$ have normalized volume equal to $1$. Then Corollary~\ref{cor:nonvertices} and Theorems~\ref{thm:special-reduced} and~\ref{thm:special-radical} imply
\[
\init_{L(\rho)}(I_{L(\mathcal A_G)}) = \bigcap_{\substack{T \in \mathcal T_\rho\\T\textup{ maximal}}} \widehat H_T, \qquad \widehat H_T \colonequals \langle x_e, y_f \colon e\not\in T\cup\Ext_\rho(T), f \in \Ext_\rho(T)\rangle.
\]
The variety $\widehat Y_T$ defined by $\widehat H_T$ is a product of projective linear subspaces
\[
\widehat Y_T = \prod_{v\in V(G)} L_{T,v}  \subseteq \prod_{v \in V(G)} \PP_v^{d_v - 1}, \qquad L_{T,v} \colonequals \left\{[\mathbf w]\in\PP_v^{d_v-1}\colon\!\! \begin{array}{l} x_e = 0 \textup{ for } s(e) = v,\ e\not\in T\cup\Ext_\rho(T),\\ y_f = 0 \textup{ for } t(f) = v,\ f\in\Ext_\rho(T)\end{array}\right\}.
\]
Thus $[\widehat Y_T] = \mathbf z(T)$, and the claim follows.
\end{proof}

We write $\mathfrak S_G(\mathbf z)$ for the multidegree of $\widehat X_G$, i.e., for the polynomial in the variables $\{z_v \colon v \in V(G)\}$ given by the formula
\[
\mathfrak S_G(\mathbf z) \overset{\rm def}= \sum_T \mathbf z(T), \qquad \mathbf z(T) = \left(\prod_{e\not\in T\cup\Ext_\rho(T)}z_{s(e)}\right)\left(\prod_{e\in\Ext_\rho(T)}z_{t(e)}\right).
\]

Our next goal is to relate $\mathfrak S_G(\mathbf z)$ to the trimmed zonotope. Recall that the \demph{integer point transform} of a polytope $P\subset \RR^n$ is the polynomial 
\[
\sum_{\mathbf p \in P\cap \ZZ^n} \mathbf x^{\mathbf p} \in \ZZ[x_1, \dots, x_n].
\]
We use the following formula for the integer point transform of (a translate of) the trimmed graphical zonotope.
\begin{thm}[{\cite[Theorem 1.4]{ms26}, cf.\ also \cite[Theorem 6.5]{kmp25}}]
\label{thm:ms26}
Let $G$ be any directed graph, and let $\indeg(G) \in \ZZ^{V(G)}$ denote the in-degree vector of $G$. Then, the integer point transform of $Z_G^- + \mathbf 1_{V(G)} + \indeg(G)$ equals
\[
P_Z(z_1, \dots, z_n) \colonequals \sum_T \mathbf x(T), \qquad \mathbf x(T) \colonequals \left(\prod_{e\in E(T)}z_{s(e)}z_{t(e)}\right)\left(\prod_{e\in \Ext_\rho(T)} z_{s(e)}\right)\left(\prod_{e\not\in E(T)\cup\Ext_\rho(T)}z_{t(e)}\right).
\]
\end{thm}
(We remark to the reader that it is not immediately clear that $\mathbf x(\cdot)$ separates trees, but it is implied by Theorem~\ref{thm:ms26}.)
\newtheorem*{thm:S_G}{Theorem~\ref{thm:S_G}}
\begin{thm:S_G}
Let $G$ be any connected directed graph. Then $\mathfrak S_G(\mathbf z)$ is the reciprocal of the integer point transform of a suitable translate of $Z_G^-$. More precisely,
\[
\mathfrak S_G(\mathbf z) = \mathbf z^\alpha \cdot \left(\sum_{\mathbf p \in Z_G^- \cap \ZZ^{V(G)}} \mathbf z^{-\mathbf p}\right),
\]
where $\alpha \colonequals \outdeg(G) - \mathbf 1_{V(G)} \in \ZZ^{V(G)}$ is a shift of the outdegree vector of $G$.
\end{thm:S_G}
\begin{proof}
Compare the formulas in Theorems~\ref{thm:multidegree} and~\ref{thm:ms26}: the monomials $\mathbf z(T)$ and $\mathbf x(T)$ satisfy
\[
\mathbf z(T) \cdot \mathbf x(T) = \mathbf z^{\deg(G)}
\]
for all $T$. It follows that
\[
\mathfrak S_G(\mathbf z) = \mathbf z^{\deg(G)} \sum_{\mathbf p \in P \cap \ZZ^{V(G)}} \mathbf z^{-\mathbf p},
\]
where $P$ is the polytope $Z_G^- + \mathbf 1_{V(G)} + \indeg(G)$. Translating $P$ backwards by $\mathbf 1_{V(G)} + \indeg(G)$ gives the claim.
\end{proof}

\section{Twisted $K$-polynomial of $\widehat X_G$}
\label{sec:G_G}
In this section, we prove Theorems~\ref{thm:G_G-top} and~\ref{thm:G_G-extra}. Theorem~\ref{thm:G_G-top} follows from Knutson--Miller--Yong's Theorem~\ref{thm:kmy} expressing (twisted) $\mathcal K$-polynomials of Stanley--Reisner rings in terms of interior faces of the Stanley--Reisner complex along with our analysis of lattice points in $Z_G$. Theorem~\ref{thm:G_G-extra} is an easy consequence of Castillo--Cid-Ruiz--Mohammadi--Monta\~no's Theorem~\ref{thm:generalized-polymatroid-support} and Knutson's Theorem~\ref{thm:mobius}.

Recall that the ring $R\colonequals \CC[\{x_e, y_e\}_{e\in E(G)}]$ is multigraded with respect to $\mdeg(x_e) = \mathbf 1_{s(e)}$ and $\mdeg(y_e) = \mathbf 1_{t(e)}$, and that $I_{L(\mathcal A_G)}$ is homogeneous with respect to $\mdeg$. Let $\mathfrak G_G(\mathbf z)$ denote the twisted $K$-polynomial of $\widehat X_G$, i.e.,
\[
\mathfrak G_G(\mathbf z) \colonequals \mathcal K(R/I_{L(\mathcal A_G)}; \mathbf 1 - \mathbf z).
\]
By Theorem~\ref{thm:kmy}, the twisted $K$-polynomial $\mathfrak G_G$ can be written as a sum over interior faces of the triangulation $\mathcal T_\rho$ of $L(\mathcal A_G)$. In our setup, this formula takes a particularly simple form. Given a set $S\subset L(\mathcal A_G)$ of vectors in the Lawrence configuration, write
\[
\mathbf z(S) \colonequals \left(\prod_{e \colon b_e \not\in S} z_{s(e)}\right)\left(\prod_{e\colon c_e \not \in S}z_{t(e)}\right).
\]
This notation extends the notation in Theorem~\ref{thm:multidegree} in the sense that if $S$ is the set of vertices of the simplex of $\mathcal T_\rho$ corresponding to a spanning tree $T$ of $G$, then $\mathbf z(S) = \mathbf z(T)$.

\begin{prop}
\label{prop:interior-faces}
We have
\[
\mathfrak G_G(\mathbf z) = \sum_{S\subset L(\mathcal A_G)} (-1)^{c - |S|}\mathbf z(S), \qquad c \colonequals |E(G)| + |V(G)| - 1,
\]
where the sum runs over the sets $S$ so that $\mathrm{conv}(S)$ forms an interior simplex of $\mathcal T_\rho$.
\end{prop}
\begin{proof}
Because the arrangement $\mathcal A_G$ is unimodular, Theorem~\ref{thm:special-radical} implies that $\init_{L(\rho)}(I_{L(\mathcal A_G)})$ is radical, i.e., $\sqrt{\init_{L(\rho)}(I_{L(\mathcal A_G)})} = \init_{L(\rho)}(I_{L(\mathcal A_G)})$. This initial ideal is a Stanley--Reisner ideal (cf.\ Remark~\ref{rem:init-is-SR}), and because $\mdeg(x_e) = \mathbf 1_{s(e)}$ and $\mdeg(y_e) = \mathbf 1_{t(e)}$ we may use Theorem~\ref{thm:kmy} to compute
\[
\mathcal K(R/I_{L(\mathcal A_G)}; \mathbf z) = \sum_S (-1)^{\codim S} \left(\prod_{e \colon b_e \not\in S} (1 - z_{s(e)})\right)\left(\prod_{e\colon c_e \not\in S}(1- z_{t(e)})\right),
\]
where the sum runs over interior faces $S$ of $\mathcal T_\rho$. Since $\dim(\mathrm{conv}(L(\mathcal A_G))) = |E(G)| + |V(G)| - 2$ and $\dim(\mathrm{conv}(S)) = |S| - 1$, substituting $z_i \mapsto 1 - z_i$ yields the claim.
\end{proof}

In order to study $\mathfrak G_G^{\mathrm{top}}(\mathbf z)$, we study simplices in the triangulation $\mathcal T_\rho$ corresponding to interior lattice points in Li--Postnikov's fine mixed subdivision of $Z_G$ into tiles $\mathcal E_{T,\rho}$.

\begin{lem}
\label{lem:zonotope-minkowski}
Let $G$ be a graph and let $T$ be a spanning tree. Every lattice point of the translated parallelepiped
\[
\mathcal E_{T,\rho} = \sum_{e\in E(T)} [0,\mathbf 1_{s(e)} - \mathbf 1_{t(e)}] + \sum_{e \in \Ext_\rho(T)}(\mathbf 1_{s(e)} - \mathbf 1_{t(e)})
\]
is of the form
\[
\left(\sum_{e\in E(T)} \alpha_e(\mathbf 1_{s(e)} - \mathbf 1_{t(e)})\right) + \sum_{e \in \Ext_\rho(T)}(\mathbf 1_{s(e)} - \mathbf 1_{t(e)}), \qquad \alpha_e \in \{0,1\}.
\]
Furthermore, for every lattice point $p \in Z_G \cap \ZZ^{V(G)}$ there exists a tree $T$ so that
\begin{equation}
\label{eqn:p's-tree}
p = \sum_{e\in S'} (\mathbf 1_{s(e)} - \mathbf 1_{t(e)}) + \sum_{e\in \Ext_\rho(T)}(\mathbf 1_{s(e)} - \mathbf 1_{t(e)}) \qquad\textup{ for some } S'\subseteq E(T).
\end{equation}
\end{lem}
\begin{proof}
Each line segment $[0,\mathbf 1_{s(e)} - \mathbf 1_{t(e)}]$ is a generalized permutahedron, so the first claim is a restatement of Lemma~\ref{lem:gp-minkowski}. Equation~\eqref{eqn:p's-tree} follows from the fact that $p \in \mathcal E_{T,\rho}$ for some $T$.
\end{proof}

\begin{lem}
\label{lem:ZG-simplex}
Let $p$ be an interior lattice point of $Z_G$, and let $T$ be a tree so that~\eqref{eqn:p's-tree} holds. Then, there is a simplex $\sigma_p \colonequals \mathrm{conv}(S_p)$ of $\mathcal T_\rho$ with vertices
\[
S_p \colonequals \{b_e \colon \mathbf 1_{s(e)} - \mathbf 1_{t(e)} \textup{ appears in~\eqref{eqn:p's-tree}}\} \sqcup \{c_e \colon \mathbf 1_{s(e)} - \mathbf 1_{t(e)} \textup{ does not appear in~\eqref{eqn:p's-tree}}\}
\]
and $\sigma_p$ intersects $Z_G$ transversely at the point $p$. In particular, $\sigma_p$ is a minimal interior face of $\mathcal T_\rho$.
\end{lem}
\begin{proof}
The formula~\eqref{eqn:p's-tree} expresses $p$ as a cell $C_p \subset \mathcal E_{T,\rho}$ in Li--Postnikov's fine mixed subdivision of $Z_G$, and under the Cayley trick the cell $C_p$ corresponds to the simplex $\sigma_p \subset T_{T,\rho}$ in the triangulation $\mathcal T_\rho$ of the Lawrence polytope $\mathrm{conv}(L(\mathcal A))$.

The simplex $\sigma_p$ is interior in the Lawrence polytope because it contains the point $p$ in the interior of $Z_G$, and it is minimal by Lemma~\ref{lem:interior-necessary}.
\end{proof}

Intuitively, Proposition~\ref{prop:interior-faces} (combined with Remark~\ref{rem:minimal-simplicies}) suggests that $\mathfrak G_G^{\mathrm{top}}(\mathbf z)$ should be a sum over lattice points in the interior of $Z_G$. It turns out, however, that $Z_G$ may have no interior lattice points at all, and the zonotopal tiling of Theorem~\ref{thm:tiling} may not have interior faces of dimension zero:
\begin{example}
Let $G$ be a graph with a leaf vertex $\ell$ (incident to an edge $e$). Observe that
\[
Z_G = Z_{G\setminus e} + [0,\mathbf 1_{s(e)} - \mathbf 1_{t(e)}]
\]
and that $\mathbf 1_{s(e)} - \mathbf 1_{t(e)}$ does not lie in the linear span of $Z_{G\setminus e}$. Thus $Z_G = Z_{G\setminus e} \times [0,\mathbf 1_{s(e)} - \mathbf 1_{t(e)}]$. Furthermore, every spanning tree $T$ of $G$ uses the edge $e$, so every face in the tiling of Theorem~\ref{thm:tiling} contains a line segment in the direction $\mathbf 1_{s(e)} - \mathbf 1_{t(e)}$.
\end{example}

It is well-known in the literature that the number of integer points in the interior of $Z_G$ is a specialization of the Tutte polynomial (e.g., using Stanley's theorem \cite[Theorem 2.2]{stanley91} on Ehrhart polynomials of zonotopes and applying Ehrhart--Macdonald reciprocity). The earliest reference we are aware of comes from the work of Holtz--Ron \cite{hr11}, Ardila--Postnikov \cite{ap10}, and Lenz \cite{lenz14} on zonotopal algebras.

Recall that the \demph{Tutte polynomial} $T_G(x,y)$ is the two-variable polynomial defined recursively via deletion-contraction by:
\begin{itemize}
\item $T_G(x,y) = T_{G\setminus e}(x,y) + T_{G/e}(x,y)$ if $e$ is not a loop or a bridge,
\item $T_G(x,y) = xT_{G\setminus e}(x,y)$ if $e$ is a bridge
\item $T_G(x,y) = yT_{G\setminus e}(x,y)$ if $e$ is a loop.
\item $T_G(x,y) = 1$ if $G$ has no edges.
\end{itemize}

\begin{lem}[{\cite[Remark 4]{lenz14}}]
\label{lem:integer-points-in-int-ZG}
The number of interior lattice points of $Z_G$ equals the Tutte polynomial specialization $T_G(0,1)$.
\end{lem}

(Lemma~\ref{lem:integer-points-in-int-ZG} is also implicit in \cite{ap10}.)

\begin{prop}
\label{prop:bridgeless-has-points}
The Tutte specialization $T_G(0,1)$ is strictly positive if and only if $G$ is bridgeless.
\end{prop}
\begin{proof}
If $G$ has a bridge $e$, then $T_G(0,1) = 0 \cdot T_{G\setminus e}(0,1) = 0$.

Conversely, observe first that $T_G(0,1) \geq 0$ for all graphs $G$ by induction on the number of edges of $G$. Now, if $G$ is bridgeless, then so is $G/e$ for all edges $e\in E(G)$, and so
\[
T_G(0,1) = T_{G\setminus e}(0,1) + T_{G/e}(0,1)
\]
is a sum of two numbers, of which the second is strictly positive by induction.
\end{proof}

We are now able to prove Theorem~\ref{thm:G_G-top}.
\newtheorem*{thm:G_G-top}{Theorem~\ref{thm:G_G-top}}
\begin{thm:G_G-top}
Let $G$ be a connected directed graph with no bridges. Then
\[
\mathfrak G_G^{\mathrm{top}}(\mathbf z) = (-1)^{|V(G)| - 1} \mathbf z^{\outdeg(G)} \cdot \left(\sum_{\mathbf p \in \mathrm{relint}(Z_G)\cap \ZZ^{V(G)}} \mathbf z^{-\mathbf p}\right),
\]
where $\mathrm{relint}(Z_G)$ denotes the relative interior of the zonotope $Z_G$.
\end{thm:G_G-top}
\begin{proof}
The monomials $\mathbf z(S)$ of maximal degree contributing to the formula for $\mathfrak G_G(\mathbf z)$ in Proposition~\ref{prop:interior-faces} are those interior simplices with the minimal number of vertices among all interior faces of $\mathcal T_\rho$. As $G$ has no bridges, Proposition~\ref{prop:bridgeless-has-points} and Lemma~\ref{lem:integer-points-in-int-ZG} implies that $Z_G$ has interior lattice points. Consequently, Lemma~\ref{lem:ZG-simplex} implies that the interior simplices of minimal dimension $S$ have exactly $|E(G)|$ many vertices.

Conversely, Lemma~\ref{lem:interior-necessary} implies that all interior simplices with exactly $|E(G)|$ many vertices are of the form $\sigma = \mathrm{conv}(\{b_e \colon e \in S\} \sqcup \{c_e \colon e \not\in S\})$ for some subset $S\subseteq E(G)$, and $\sigma \cap Z_G$ is equal to the lattice point $\sum_{e\in S} (\mathbf 1_{s(e)} - \mathbf 1_{t(e)})$. Thus, an interior simplex with exactly $|E(G)|$ many vertices are precisely those equal to $\sigma_p$ for $p$ an interior lattice point of $Z_G$.

It remains to compute $\mathbf z(S_p)$. To this end, write~\eqref{eqn:p's-tree} in the form
\[
p = \sum_{e\in S} (\mathbf 1_{s(e)} - \mathbf 1_{t(e)}), \qquad \textup{ for some } S\subseteq E(G)
\]
and note that
\begin{align*}
\outdeg(G) - p &= \sum_{e\in E(G)} \mathbf 1_{s(e)} - p \\&= \sum_{e\not\in S}\mathbf 1_{s(e)} + \sum_{e\in S}\mathbf 1_{t(e)}.
\end{align*}
Lemma~\ref{lem:ZG-simplex} asserts that $e\not\in S$ if and only if $b_e\not\in S_p$ and that $e\in S$ if and only if $c_e\not\in S_p$, so in particular
\begin{align*}
\mathbf z(S_p) &= \left(\prod_{e \colon b_e \not\in S_p} z_{s(e)}\right)\left(\prod_{e\colon c_e \not \in S_p}z_{t(e)}\right) \\&= \mathbf z^{\outdeg(G) - p}.
\end{align*}
It follows that
\begin{align*}
\mathfrak G_G^{\mathrm{top}} &= \sum_{\mathbf p \in \mathrm{relint}(Z_G)\cap \ZZ^{V(G)}} (-1)^{c - |E(G)|} \cdot \mathbf z(S_{\mathbf p})\\&= (-1)^{|V(G)| - 1}\mathbf z^{\outdeg(G)} \cdot \left(\sum_{\mathbf p \in \mathrm{relint}(Z_G)\cap \ZZ^{V(G)}} \mathbf z^{-\mathbf p}\right).\qedhere
\end{align*}
\end{proof}

\subsection*{Multiplicity-freeness and twisted $K$-polynomials}

In this subsection we use multiplicity-freeness of $\widehat X_G$ to study $\mathfrak G_G(\mathbf z)$. More specifically, we apply a theorem of Castillo--Cid-Ruiz--Mohammadi--Monta\~no (\cite[Theorem A]{ccmm22}) to obtain that the homogenization $\widetilde{\mathfrak G}_G(\mathbf z)$ has M-convex support, and a formula of Knutson \cite[Theorem 3]{knutson09} (see also~\cite{ps24}) to obtain $\mathfrak G_G(\mathbf z)$ from M\"obius inversion (Corollary~\ref{cor:mobius-explicit}). Theorem~\ref{thm:G_G-extra} is essentially an immediate consequence of these two theorems.

\begin{prop}
\label{prop:tilde-M-convex}
The support of the homogenization $\widetilde{\mathfrak G}_G(\mathbf z)$ is M-convex, and in particular $\supp(\widetilde{\mathfrak G}_G(\mathbf z))$ is saturated.
\end{prop}
\begin{proof}
Theorem~\ref{thm:S_G} implies that $\widehat X_G \subset\prod_{v\in V(G)} \PP^{\deg(v)-1}$ is multiplicity-free. Thus Theorem~\ref{thm:generalized-polymatroid-support} implies that $\mathfrak G_G(\mathbf z)$ has generalized polymatroid support, and that the support is in particular saturated.
\end{proof}

\begin{prop}
\label{prop:G_G-mobius}
The coefficients of $\mathfrak G_G(\mathbf z)$ can be computed via M\"obius inversion:
\[
\mathfrak G_G(\mathbf z) = \sum_{\mathbf q\in P_G}\mu(\mathbf q) \mathbf z^{\mathbf q},
\]
where $P_G$ is the poset $\{\mathbf q \in \ZZ_{\geq 0}^{V(G)}\colon \mathbf q \leq \deg(G) - \mathbf 1 \textup{ and } \mathbf z^{\mathbf q}\textup{ is divisible by some monomial in } \mathfrak S_G(\mathbf z)\}$ and $\mu$ is the M\"obius function of $P_G$.
\end{prop}
\begin{proof}
Theorem~\ref{thm:S_G} implies that $\widehat X_G \subset \prod_{v\in V(G)} \PP^{\deg(v) - 1}$ is multiplicity-free. In our setting, \cite[Theorem 3]{knutson09} (Corollary~\ref{cor:mobius-explicit}) states that
\[
\mathfrak G_G(\mathbf z) = \sum_{\mathbf q \in P_G} \mu(\mathbf q) \mathbf z^{\mathbf q},
\]
where $\mu$ is the M\"obius function and $P_G$ is the poset of exponent vectors dominating a vector in the support of $\mathfrak S_G(\mathbf z)$ and dominated by the vector $\deg(G) - \mathbf 1$.
\end{proof}

\newtheorem*{thm:G_G-extra}{Theorem~\ref{thm:G_G-extra}}
\begin{thm:G_G-extra}
The homogenization $\widetilde{\mathfrak G}_G(\mathbf z)$ of the twisted $K$-polynomial of $\widehat X_G$ has M-convex support, and in particular satisfies the saturated Newton polytope property. Furthermore, the coefficients of $\mathfrak G_G(\mathbf z)$ can be computed via M\"obius inversion:
\[
\mathfrak G_G(\mathbf z) = \sum_{\mathbf q\in P_G}\mu(\mathbf q) \mathbf z^{\mathbf q},
\]
where $P_G$ is the poset $\{\mathbf q \in \ZZ_{\geq 0}^{V(G)}\colon \mathbf q \leq \deg(G) - \mathbf 1\textup{ and } \mathbf z^{\mathbf q}\textup{ is divisible by some monomial in } \mathfrak S_G(\mathbf z)\}$ and $\mu$ is the M\"obius function of $P_G$.
\end{thm:G_G-extra}
\begin{proof}
Combine Propositions~\ref{prop:tilde-M-convex} and~\ref{prop:G_G-mobius}.
\end{proof}

\section{Flat arrangements and log-concavity}
\label{sec:log-concave}
The goal of this section is to prove Theorem~\ref{thm:log-concave}. \textbf{For this section}, we assume that $\mathcal A$ is \demph{flat}, meaning that there exists a linear functional $h\colon \RR^n\to\RR$ with $h(a_i) = 1$ for all $i$. We also assume that $\mathcal A$ is full-dimensional. Though we require flatness in this section (and not in Section~\ref{sec:graphical}), the arrangement need not be unimodular (in contrast to Section~\ref{sec:graphical}).

To start, let us assume $\mathcal A = \{a_1, \dots, a_N\}$ is \demph{rational}, meaning that $a_i \in \QQ^n$ for all $i$, and let $\Lambda_{\mathcal A} \colonequals \sum_{i=1}^N\ZZ a_i$ denote the lattice generated by $\mathcal A$. We first prove Theorem~\ref{thm:log-concave} for rational arrangements, and then treat the general case by an approximation argument.

\begin{lem}
\label{lem:bideg}
The ideal $I_{L(\mathcal A)}$ is homogeneous with respect to the bigrading defined by $\bideg(x_i)\colonequals (1,0)$ and $\bideg(y_i)\colonequals (0,1)$.
\end{lem}
\begin{proof}
The ideal $I_{L(\mathcal A)}$ is generated by binomials of the form
\begin{equation*}
\label{eqn:bideg}
\tag{$\dagger$}
\mathbf x^\alpha \mathbf y^{\alpha'} - \mathbf x^\beta \mathbf y^{\beta'}, \qquad \textup{ where }\sum_{i=1}^N \alpha_i b_i + \sum_{i=1}^N \alpha_i' c_i = \sum_{i=1}^N \beta_i b_i + \sum_{i=1}^N \beta_i'c_i.
\end{equation*}
Extend $h\colon \RR^n \to \RR$ to a linear functional $\tilde h \colon \RR^n\times\RR^N \to \RR$ by setting $\tilde h(a,z) \colonequals h(a)$. Note in particular that $\tilde h(b_i) = 1$ and $\tilde h(c_i) = 0$.

Applying $\tilde h$ to the condition on $\alpha,\alpha',\beta,\beta'$ in~\eqref{eqn:bideg} yields $|\alpha| = |\beta|$. Then, applying the linear functional $(a,z)\mapsto |z|$ to the condition on $\alpha,\alpha',\beta,\beta'$ in~\eqref{eqn:bideg} yields $|\alpha| + |\alpha'| = |\beta| + |\beta'|$. It follows that $(|\alpha|,|\alpha'|) = (|\beta|,|\beta'|)$, and so $I_{L(\mathcal A)}$ is generated by homogeneous polynomials with respect to $\bideg$.
\end{proof}

\begin{rem}
If $G$ is a bipartite graph with vertex bipartition $V(G) = V(G)_L \sqcup V(G)_R$, and all edges oriented from $V(G)_L$ to $V(G)_R$, then $\mdeg$ is a refinement of $\bideg$, by taking $\mathbf 1_v \mapsto (1,0)$ if $v \in V(G)_L$ and $\mathbf 1_v\mapsto (0,1)$ if $v\in V(G)_R$. It will follow that, in this case, the homogenization of $f_{\mathcal A}(t)$ is a specialization of $\mathfrak S_G(\mathbf z)$.
\end{rem}

\begin{defn}
Let $X_{\mathcal A}\subseteq \PP_x^{N-1}\times\PP_y^{N-1}$ denote the biprojective variety associated to $I_{L(\mathcal A)}$.
\end{defn}

\begin{prop}
\label{prop:XA-formula}
The Chow class $[X_{\mathcal A}] \in A^*(\PP^{N-1}_x \times \PP^{N-1}_y)$ is equal to
\[
[X_{\mathcal A}] = \sum_B \Vol_{\mathrm{norm}}(T_{B,\rho})[H_B], \qquad H_B = \left\{([\mathbf x], [\mathbf y]) \in \PP_x^{N-1}\times\PP_y^{N-1} \colon\!\! \begin{array}{l} x_i = 0 \textup{ for } i \not\in B\cup \Ext_\rho(B),\\ y_j = 0 \textup{ for } j \in \Ext_\rho(B)\end{array}\right\}.
\]
\end{prop}
\begin{proof}
By Lemma~\ref {lem:init-is-monomial}, the initial ideal is monomial, and Proposition~\ref{prop:omega} implies that there is a term order $\prec_\omega$ so that $\init_{L(\rho)}(I_{L(\mathcal A)}) = \init_{\prec_\omega}(I_{L(\mathcal A)})$. Hence, there is a Gr\"obner degeneration
\[
X_{\mathcal A} \rightsquigarrow X_0
\]
where $X_0$ is the projective scheme associated to $\init_{L(\rho)}(I_{L(\mathcal A)})$. Theorems~\ref{thm:special-reduced} and~\ref{thm:special-multiplicity} imply that
\[
[X_0] = \sum_B \Vol_{\mathrm{norm}}(T_{B,\rho})[H_B]
\]
because the nonvertices of $T_{B,\rho}$ are $b_i$ for $i \not\in B\cup \Ext_\rho(B)$ and $c_j$ for $j\in\Ext_\rho(B)$. The claim follows from the fact that Chow class is preserved under Gr\"obner degeneration.
\end{proof}

Let us denote by $\widetilde f_{\mathcal A}(q,t)$ the homogenization of $f_{\mathcal A}(t)$. Concretely,
\[
\widetilde f_{\mathcal A}(q,t) = \sum_{B\in\mathcal B(\mathcal A)}\Vol(\Pi_B) q^{N-n-\ext_\rho(B)} t^{\ext_\rho(B)}.
\]
\begin{thm}
\label{thm:main-geometric}
Let $q$ and $t$ denote the hyperplane classes of $\PP^{N-1}_x$ and $\PP^{N-1}_y$ respectively. Then
\[
[X_{\mathcal A}] = \frac1{\covol(\Lambda_{\mathcal A})}\cdot \widetilde f_{\mathcal A}(q,t).
\]
\end{thm}
\begin{proof}
Proposition~\ref{prop:XA-formula} and Lemma~\ref{lem:volumes} gives
\[
[X_{\mathcal A}] = \frac1{\covol(\Lambda_{\mathcal A})}\cdot \sum_B \Vol(\Pi_B)[H_B].
\]
The variety $H_B$ is a product of linear subspaces of $x$-codimension $N-n-\ext_\rho(B)$ and of $y$-codimension $\ext_\rho(B)$. Thus $[H_B] = q^{N-n-\ext_\rho(B)} t^{\ext_\rho(B)}$, and the claim follows.
\end{proof}
\begin{cor}
\label{cor:log-concavity-rational}
Let $\mathcal A$ be a rational, full-dimensional, flat vector arrangement. Then the coefficients of $f_{\mathcal A}(t)$ are log-concave.
\end{cor}
\begin{proof}
The ideal $I_{L(\mathcal A)}$ is prime, so the variety $X_{\mathcal A}$ is irreducible, so \cite[Theorem 21]{huh12} implies that the coefficients of $\widetilde f_{\mathcal A}(q,t)$ are log-concave. (De-homogenizing, $f_{\mathcal A}(t) = \widetilde f_{\mathcal A}(1,t)$ is then log-concave.)
\end{proof}

Strictly speaking, the proof of Corollary~\ref{cor:log-concavity-rational} requires only one direction of the if-and-only-if statement in \cite[Theorem 21]{huh12}. The direction we need is a consequence of the Khovanskii--Teissier inequalities \cite[Example 1.6.4]{lazarsfield04}.

\subsection*{General arrangements}
We remove the rationality assumption on $\mathcal A$ in Corollary~\ref{cor:log-concavity-rational}, and prove Theorem~\ref{thm:log-concave}. We need to approximate $\mathcal A$ by rational arrangements $\mathcal A^{(m)}$ and argue that $f_{\mathcal A^{(m)}}(t) \to f_{\mathcal A}(t)$.

\begin{prop}
\label{prop:approx}
Let $\mathcal A$ be a full-dimensional flat vector arrangement in $\RR^n$. Then there is a sequence of rational, full-dimensional, flat vector arrangements $\mathcal A^{(m)}$ such that $f_{\mathcal A^{(m)}}(t) \to f_{\mathcal A}(t)$ coefficientwise.
\end{prop}
\begin{proof}
The polynomial $f_{\mathcal A}(t)$ is invariant, up to scaling, under simultaneous action by $\mathrm{GL}_n(\RR)$, so we may assume $h\colon \RR^n\to\RR$ is given by $(x_1, \dots, x_n)\mapsto x_n$, and hence that $a_i$ takes the form $a_i = (p_i, 1)$ for $p_i \in \RR^{n-1}$.

Choose a sequence of vectors $a_i^{(m)}$ of the form $a_i^{(m)} = (p_i^{(m)},1)$ for $p_i^{(m)} \in \QQ^{n-1}$ with $p_i^{(m)}\to p_i$. Given any subset $S$ of $\mathcal A$ of size $n$, write $\Pi_S^{(m)}$ for the zonotope spanned by $\{a_i^{(m)}\colon i \in S\}$ and observe that $\Vol(\Pi_S^{(m)}) \to \Vol(\Pi_S)$. Hence, by taking $S$ to be any basis of $\mathcal A$, the rational and flat arrangement $\mathcal A^{(m)}$ is also full-dimensional for $m$ sufficiently large; even more, a basis $B$ of $\mathcal A$ is also a basis for $\mathcal A^{(m)}$ for all $m$ sufficiently large.

Choose a generic framing $\rho\in\RR^N$ for $\mathcal A$ so that the coordinates of $\rho$ are linearly independent over $\QQ$. Then $\rho$ is automatically generic for the rational arrangements $\mathcal A^{(m)}$, as circuit vectors in $\mathcal A^{(m)}$ can be chosen so that all of their components are rational.

We next argue that for a fixed basis $B$ of $\mathcal A$, any fixed $j\not\in B$ is externally semiactive in $\mathcal A$ if and only if it is externally semiactive in $\mathcal A^{(m)}$ for all $m\gg0$. To this end, let $C$ and $C^{(m)}$ be the fundamental circuits of $j$ with respect to $B$ in $\mathcal A$ and $\mathcal A^{(m)}$ respectively. View the normalized circuit vectors $\lambda^{(C,j)}$ and $\lambda^{(C^{(m)},j)}$ as vectors in the kernels of
\[
\begin{aligned}
\varphi\colon \RR^{B\cup\{j\}} &\twoheadrightarrow \RR^n\\[-3pt]
\mathbf 1_i&\mapsto a_i \vphantom{a_i^{(m)}}
\end{aligned}
\qquad\textup{ and } \qquad
\begin{aligned}
\varphi^{(m)}\colon \RR^{B\cup\{j\}} &\twoheadrightarrow \RR^n\\[-3pt]
\mathbf 1_i&\mapsto a_i^{(m)}
\end{aligned}
\]
respectively. Since $\varphi^{(m)} \to \varphi$ coordinatewise, we have $\lambda^{(C^{(m)},j)}\to \lambda^{(C,j)}$. Thus
\[
\langle \rho, \lambda^{(C^{(m)},j)}\rangle \to \langle \rho, \lambda^{(C,j)}\rangle
\]
and in particular $\langle \rho, \lambda^{(C^{(m)},j)}\rangle$ will eventually have the same sign as $\langle \rho, \lambda^{(C,j)}\rangle$. We conclude that $j$ is externally semiactive in $\mathcal A$ if and only if it is externally semiactive in $\mathcal A^{(m)}$ for all $m\gg0$.

Fix $\varepsilon > 0$ and take a sufficiently large $M$ so that, whenever $m \geq M$, 
\begin{itemize}
\item Every basis of $\mathcal A$ is also a basis of $\mathcal A^{(m)}$,
\item For every basis $B$ and every $j\not\in B$, the element $j$ is externally semiactive in $\mathcal A$ if and only if $j$ is externally semiactive in $\mathcal A^{(m)}$, and
\item $|\Vol(\Pi_S^{(m)}) - \Vol(\Pi_S)|\leq \varepsilon$ for all subsets $S$ of size $n$.
\end{itemize}
Then, whenever $m\geq M$, we have
\[
\left|[t^k] f_{\mathcal A^{(m)}}(t) - [t^k] f_{\mathcal A}(t)\right| \leq \sum_{\substack{S\subseteq[N]\\|S| = n}}|\Vol(\Pi_S^{(m)}) - \Vol(\Pi_S)| \leq \varepsilon \cdot \binom Nn.\qedhere
\]
\end{proof}
\newtheorem*{thm:log-concave}{Theorem~\ref{thm:log-concave}}
\begin{thm:log-concave}
Let $\mathcal A$ be a full-dimensional flat vector arrangement. Then the coefficients of $f_{\mathcal A}(t)$ are log-concave.
\end{thm:log-concave}
\begin{proof}
By Proposition~\ref{prop:approx}, there exist rational arrangements $\mathcal A^{(m)}$ so that $f_{\mathcal A^{(m)}}(t) \to f_{\mathcal A}(t)$. By Corollary~\ref{cor:log-concavity-rational}, each $f_{\mathcal A^{(m)}}(t)$ has log-concave coefficients. Because log-concavity is a closed condition, $f_{\mathcal A}(t)$ has log-concave coefficients too.
\end{proof}

\end{document}